\documentclass[final,3p,times]{elsarticle}

\usepackage[hidelinks]{hyperref}
\usepackage{amsmath,amssymb,amsthm}
\usepackage{epsfig}
\usepackage{float}
\usepackage{caption}
\usepackage{subcaption}
\usepackage{moreverb}
\usepackage{epstopdf}
\usepackage{booktabs}

\usepackage{graphicx}

\usepackage[T1]{fontenc}
\usepackage[utf8]{inputenc}

\newtheorem{theorem}{Theorem}[section]

\newtheorem{definition}[theorem]{Definition}

\theoremstyle{remark}
\newtheorem{remark}[theorem]{Remark}

\usepackage[dvipsnames]{xcolor}
\newcommand\revOne[1]{\textcolor{black}{#1}}
\newcommand\revTwo[1]{\textcolor{black}{#1}}
\newcommand\revThree[1]{\textcolor{black}{#1}}

\journal{Journal of Computational Physics}

\newcommand{\orcid}[1]{\href{https://orcid.org/#1}{\includegraphics[width=10pt]{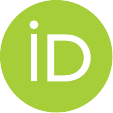}}}

\begin{document}

\begin{frontmatter}

\title{\revOne{Boundary-optimized closures for diagonal-norm upwind SBP operators}}

\author[upp]{Ken Mattsson\orcid{0000-0003-0843-3321}}

\author[upp]{David Niemel\"{a}} 

\author[liu]{Andrew R. Winters\corref{cor1}\orcid{0000-0002-5902-1522}}
\ead{andrew.ross.winters@liu.se}
\cortext[cor1]{Corresponding author}

\affiliation[upp]{organization={Uppsala University, Department of Information Technology},
            city={Uppsala},
            postcode={75105}, 
            country={Sweden}}

\affiliation[liu]{organization={Linköping University, Department of Mathematics},
            city={Linköping},
            postcode={58183}, 
            country={Sweden}}

\begin{abstract}
By employing non-equispaced grid points near boundaries, boundary-optimized upwind finite-difference operators of orders up to nine are developed. The boundary closures are constructed within a diagonal-norm summation-by-parts (SBP) framework, ensuring linear stability on piecewise curvilinear multiblock grids. \revTwo{For linear problems, this stability is inherited from the diagonal norm SBP framework provided the approximation satisfies the metric identities (i.e. it is free-stream preserving). For nonlinear problems, the flux-vector splitting should be linear in the metric terms to inherit these stability properties.} Boundary and interface conditions are imposed using either weak enforcement through simultaneous approximation terms (SAT) or strong enforcement via the projection method.

\revTwo{The proposed operators yield significantly improved accuracy compared with SBP operators constructed on equidistant grids.}  The resulting SBP--SAT and SBP--projection discretizations produce fully explicit systems of ordinary differential equations. The accuracy and stability properties of the proposed operators are demonstrated through numerical experiments for linear hyperbolic problems in one spatial dimension and for the compressible Euler equations in two spatial dimensions.
\end{abstract}


\end{frontmatter}

\section{Introduction}\label{sec:intro}

High-order numerical methods are well known to capture transient phenomena more efficiently than first- and second-order methods, as they allow a substantial reduction in the number of degrees of freedom required to achieve a prescribed error tolerance. In particular, high-order finite difference methods are well suited for such problems; see the pioneering work by Kreiss and Oliger \cite{KreissOliger72}. For long-time simulations, it is essential to employ finite difference approximations that do not introduce artificial growth when the underlying partial differential equation is energy bounded, i.e., methods that satisfy \textit{time stability} \cite{Gustafsson2007}.

For well-posed initial-boundary value problems (IBVPs), a robust and widely used high-order finite difference framework combines summation-by-parts (SBP) operators \cite{KreissScherer74,Strand94,MattssonNordstrom04} with either the simultaneous approximation term (SAT) method or the projection method for imposing boundary conditions. The SAT formulation can also be interpreted as a two-point numerical flux formulation; see \cite{ranocha2025robustness}. The SBP property is independent of the underlying discretization paradigm and has been successfully applied in finite difference, finite volume \cite{ham2006accurate,nordstrom2003finite}, finite element \cite{abgrall2020analysis,Hanifa2025}, and discontinuous Galerkin methods \cite{winters2021,gassner2013skew}.

For problems with variable coefficients, such as spatially varying material parameters or discretizations on curvilinear grids, diagonal-norm SBP operators are typically required to guarantee stability unless artificial dissipation (AD) is added to the scheme \cite{Svard04,MattssonAlmquist13}. Recent developments of the SBP-SAT methodology include \revTwo{\cite{VidarMattsson_RV_2021,VidarMattsson_D2_2023, ranocha2025robustness, del2019extension, chan2018discretely, ranocha2017extended}} and the review papers \cite{DelReyFernandez2014a,Svard2014}. Corresponding developments of the SBP-projection methodology are presented in \cite{MattssonOlsson2018,ErikssonMattsson22}.

Most SBP operators available in the literature (see, e.g., \cite{KreissScherer74,Strand94,MattssonNordstrom04,Mattsson11,Mattsson14,Linders1109585,MattssonAlmquistWeide2018,DelReyFernandez2014}) are based on central finite difference stencils with specially designed one-sided boundary closures that mimic the continuous integration-by-parts property in a discrete norm. Traditional diagonal-norm SBP operators, here referring to operators constructed on equidistant grids, suffer from reduced formal order of accuracy near boundaries. It has been shown that diagonal-norm SBP operators where all finite difference stencils have accuracy of at least order $p$ are associated with quadrature rules of degree at least $2p-1$\revTwo{, see, e.g., \cite{DelReyFernandez2014,hicken2013summation}}.

To mitigate this boundary accuracy reduction, boundary-optimized SBP operators for first derivatives were introduced in \cite{MattssonAlmquistCarpenter13} and later improved in \cite{MattssonAlmquistWeide2018}. These operators retain $p^{\text{th}}$-order accuracy near boundaries while significantly reducing the leading-order truncation error constant. This improvement is achieved by employing stencils constructed on non-equispaced grid points near the boundaries. More recently, boundary-optimized SBP operators for second derivatives with variable coefficients were presented in \cite{VidarMattsson_D2_2023}.

For nonlinear hyperbolic and hyperbolic-parabolic IBVPs, such as the compressible Euler and Navier--Stokes equations, the inclusion of artificial dissipation is often necessary to suppress spurious oscillations. However, the construction of high-order accurate AD operators is nontrivial and typically involves tuning free parameters. The added AD must preserve linear stability, accuracy, and non-stiffness, which requires carefully designed boundary closures. The incorporation of AD within traditional central-difference SBP-SAT discretizations was first analyzed in \cite{MattssonSvard04}. The diagonal-norm upwind SBP operators introduced in \cite{Mattsson17} naturally provide high-order artificial dissipation, effectively damping high-frequency oscillations when combined with standard flux-splitting techniques; see, e.g., \cite{Mattsson17,LundgrenMattsson_2020,duru2024,hew2025strongly}. For strongly nonlinear phenomena such as shocks, additional dissipation mechanisms are required. An extension of residual-based artificial viscosity methods to the finite difference SBP-SAT framework is presented in \cite{VidarMattsson_RV_2021}.

In addition to their dissipative properties, the upwind SBP operators introduced in \cite{Mattsson17} provide improved dispersion characteristics for wave propagation problems compared with central-difference SBP operators \cite{MattssonOssian17,williams2024full}.

The primary contribution of the present work is the derivation of diagonal-norm upwind SBP operators with boundary closures optimized for accuracy. The proposed operators are constructed using non-equispaced boundary grids and are shown to outperform both the standard upwind SBP operators introduced in \cite{Mattsson17} and the boundary-optimized central-difference SBP operators presented in \cite{MattssonAlmquistWeide2018}. \revTwo{For several operator orders, superconvergence is observed for the newly derived boundary-optimized SBP operators in addition to observed reduced errors at moderate resolution compared to other available SBP operators.}

In Section \ref{sec:FD}, the novel upwind SBP definition is presented together with the key steps of its derivation. Section \ref{sec:1-D_1} introduces the SBP-projection method for linear hyperbolic problems in one dimension and highlights the differences between traditional and upwind discretizations. In Section \ref{sec:Computations_1D}, the accuracy and stability properties of the optimized upwind SBP operators are verified and compared with previously derived upwind and central-difference SBP operators through numerical simulations of one-dimensional hyperbolic problems. Section \ref{sec:Computations_2D} presents long-time numerical simulations of an analytical two-dimensional Euler vortex and the Kelvin--Helmholtz instability on multiblock grids to further verify accuracy and stability. Finally, Section \ref{sec:Conclusions} summarizes the work. 
\revThree{The coefficients of the optimized upwind SBP operators are included online as supplementary data for this paper, as described in \ref{AppA}.}

\section{The finite difference method}\label{sec:FD}

The following definitions are required for the present study. We consider real-valued vector functions with $k$ components ${\bf u},{\bf v}\in L^2([0,1])^k$, where
\[
{\bf u}^T=[u^{(1)},u^{(2)},\ldots,u^{(k)}], \qquad
{\bf v}^T=[v^{(1)},v^{(2)},\ldots,v^{(k)}].
\]
The continuous weighted inner product on the interval $[x_l,x_r]$ is defined by
\[
({\bf u},{\bf v})=\int_{x_l}^{x_r}{\bf u}^T(x)\,{\bf C}(x)\,{\bf v}(x)\,\mathrm{d}x,
\qquad {\bf C}(x)={\bf C}^T(x)>0,
\]
and the corresponding weighted norm is $\|{\bf u}\|_C^2=({\bf u},{\bf u})$.

We discretize the domain $0\le x\le 1$ using $m$ grid points,
\[
{\bf x}=[x_1,\,x_2,\ldots,x_m]^T,
\]
where ${\bf x}$ denotes the vector of grid points. For an equidistant grid,
\[
x_i=(i-1)h,\qquad i=1,2,\ldots,m,\qquad h=\frac{1}{m-1}.
\]

To improve boundary accuracy, we allow the two grid points closest to each boundary to be non-equidistant. We introduce $d=d_1+d_2$, where $h d_1$ is the distance between the first and second grid points and $h d_2$ is the distance between the second and third grid points. The grid spacing is defined as
\[
h=\frac{1}{2d+(m-5)},
\]
where $m$ denotes the total number of grid points. The resulting non-equidistant grid is
\begin{equation}\label{eq:gridpoints}
{\bf x}=[0,\,d_1h,\,dh,\,(d+1)h,\,(d+2)h,\ldots,(d+(m-6))h,\,1-dh,\,1-d_1h,\,1]^T,
\end{equation}
where $d_1$ and $d_2$ are chosen to minimize boundary truncation errors.

For a system of equations, we define the concatenated discrete solution vector by
\[
v^T=[v^{(1)},v^{(2)},\ldots,v^{(k)}],
\]
where $v^{(j)}=[v^{(j)}_1,v^{(j)}_2,\ldots,v^{(j)}_m]$ denotes the discrete solution vector associated with the $j^{\text{th}}$ equation. 
\revOne{We make use of the Kronecker product
\[
C\otimes D=
\begin{bmatrix}
c_{0,0}D & \cdots & c_{0,q-1}D\\
\vdots & & \vdots\\
c_{p-1,0}D & \cdots & c_{p-1,q-1}D
\end{bmatrix},
\]
where $C$ is a $p\times q$ matrix and $D$ is an $m\times n$ matrix. Two useful identities are
\[
(A\otimes B)(C\otimes D)=(AC)\otimes(BD),\qquad (A\otimes B)^T=A^T\otimes B^T.
\]
}

For discrete real-valued vector functions $u,v\in\mathbb{R}^{km}$, we define the discrete inner product
\[
(u,v)_{H_kC}=u^T\,C H_k\,v,
\]
where $H_k=I_k\otimes H$, $H$ is a symmetric positive-definite $m\times m$ matrix, and $C$ denotes the projection of ${\bf C}(x)$ onto block-diagonal form. If ${\bf C}(x)$ is diagonal, then $C$ is diagonal. Here $I_k$ denotes the $k\times k$ identity matrix. The corresponding norm is $\|v\|_{H_kC}^2=v^T\,C H_k\,v$.

\begin{remark}
The matrix product $C H_k$ defines a norm if and only if it is positive definite. For variable coefficient matrices ${\bf C}(x)$, this property can only be guaranteed when the norm matrix $H$ is diagonal; see \cite{Svard04}. Variable coefficients commonly arise from spatially varying material parameters or coordinate transformations associated with curvilinear grids.
\end{remark}
The unit basis vectors $e_1$ and $e_m$, and the matrix $B$, are defined by
\begin{equation}\label{eq:Def_matrix}
e_1=[1,0,\ldots,0]^T,\qquad e_m=[0,\ldots,0,1]^T,\qquad B=e_me_m^T-e_1e_1^T.
\end{equation}

In the present study, SBP operators constructed on equidistant grids are referred to as \textit{traditional} SBP operators. Traditional central-difference first-derivative SBP operators were introduced in \cite{KreissScherer74,Strand94}, and traditional second-derivative SBP operators were introduced in \cite{Mattsson11}.

\begin{definition}\label{definition:SBP_D1}
We define a difference operator
\[
D_1=H^{-1}\left(Q+\frac{B}{2}\right)
\]
approximating $\partial/\partial x$ using a $p^{\text{th}}$-order accurate interior stencil as a $p^{\text{th}}$-order first-derivative SBP operator if $H$ is symmetric positive definite and $Q+Q^T=0$.
\end{definition}

\subsection{Upwind SBP operators}\label{sec:Upwind_SBP}

We define upwind SBP operators \cite{Mattsson17} as follows.
\begin{definition}\label{definition:SBP_Dp2}
We define
\[
D_{+}=H^{-1}\left(Q_+ +\frac{B}{2}\right),\qquad
D_{-}=H^{-1}\left(Q_-+\frac{B}{2}\right),
\]
approximating $\partial/\partial x$ using $p^{\text{th}}$-order accurate interior stencils, as $p^{\text{th}}$-order diagonal-norm upwind SBP operators if $H$ is diagonal and positive definite, $Q_+ + Q_-^T=0$, and $Q_+ + Q_+^T=2S$ is negative semi-definite.
\end{definition}

Traditional upwind SBP operators of both odd and even orders were derived in \cite{Mattsson17} up to ninth order. The following relations are frequently used:
\begin{equation}\label{eq:Relations_1D}
\begin{aligned}
\frac{D_{+}+D_{-}}{2}
&=H^{-1}\left(\frac{Q_{+}+Q_{-}}{2}+\frac{B}{2}\right)
=H^{-1}\left(Q+\frac{B}{2}\right)=D_1,\\[0.1cm]
\frac{D_{+}-D_{-}}{2}
&=H^{-1}\left(\frac{Q_{+}-Q_{-}}{2}\right)
=H^{-1}\left(\frac{Q_{+}+Q_{+}^T}{2}\right)=H^{-1}S,\\[0.1cm]
HD_{-}&=-(D_{+})^T H+B,\qquad HD_{+}=-(D_{-})^T H+B,
\end{aligned}
\end{equation}
where $Q=-Q^T$ and $S=S^T\le 0$. Hence, $(D_{+}+D_{-})/2$ yields a central SBP approximation of the first derivative, while $(D_{+}-D_{-})/2$ defines a negative semi-definite operator acting as artificial dissipation.

Next, we derive boundary-optimized upwind SBP operators.

\subsection{Boundary-optimized upwind SBP operators}\label{sec:1-D_1T}

To motivate the construction, we introduce additional notation. Let
\[
\hat{D}=H^{-1}\left(\hat{Q}+\frac{B}{2}\right)
\]
denote either a traditional first-derivative SBP operator $D_1$ (Definition~\ref{definition:SBP_D1}) or an upwind SBP operator $D_{\pm}$ (Definition~\ref{definition:SBP_Dp2}). The following two definitions form the foundation of the present work.

\begin{definition}\label{definition:accuracy}
Let ${\bf x}^q$ denote the grid function obtained by evaluating the polynomial $\frac{x^q}{q!}$ at the grid points. The $q^{\text{th}}$-order error vector is defined by
\begin{equation}\label{eq:error_vec}
{\bf e}_{(q)} = Hq{\bf x}^{\revOne{\max(0,q-1)}}-\left(\hat{Q}+\frac{B}{2}\right){\bf x}^{q}.
\end{equation}
We say that $\hat{D}$ is $p^{\text{th}}$ order accurate if ${\bf e}_{(q)}$ vanishes for \revOne{$q=0,\ldots,p$} in the interior and for \revOne{$q=0,\ldots,s/2$} at the boundaries, where $s=p$ for even $p$ and $s=p-1$ for odd $p$.
\end{definition}

\begin{definition}\label{definition:error}
The discrete $\ell_2$-norm of ${\bf e}_{(q)}$ in \eqref{eq:error_vec} is defined by
\[
\|{\bf e}_{(q)}\|_{h}^2 = h\,{\bf e}_{(q)}^T{\bf e}_{(q)},
\]
where $h$ is the grid spacing.
\end{definition}
\revOne{The quantity in Definition~\ref{definition:error} is the algebraic error measure used in the operator construction and should not be interpreted as an SBP-quadrature approximation of the continuous $L^2$ norm of the pointwise truncation error. For the nonlinear parameter optimization described below, we use an unscaled reference grid with interior spacing $h=1$ and $m=21$ grid points. Thus, the factor $h$ in Definition~\ref{definition:error} is unity during the Maple optimization; the resulting operators are subsequently scaled to the physical grid spacing.}

The formal boundary accuracy does not, by itself, fully determine the observed convergence rate. A careful error analysis (not repeated here; see \cite{SVARD2019108819}) provides guidance on the expected rates, summarized in the following remark.

\begin{remark}
Sv\"{a}rd and Nordstr\"{o}m \cite{SVARD2019108819} showed that a stable approximation of an IBVP involving derivatives up to order $q$ yields a convergence rate of order $q+b$, where $b$ is the boundary accuracy. Let $p$ denote the interior accuracy of a diagonal-norm SBP operator (traditional or upwind). \revThree{So, if $q+b > p$ then the convergence order will be $p$, otherwise the convergence rate is $q+b$ dominated by the boundary closure.} For first-order hyperbolic problems this implies a convergence rate of order $b+1$. For even $p$, diagonal-norm boundary closures are restricted to accuracy $b=p/2$ \cite{MattssonNordstrom04}, implying that ${\bf e}_{(p/2+1)}$ is the leading-order error term. For odd $p$, the leading-order error term is ${\bf e}_{((p-1)/2+1)}$. For example, the eighth- and ninth-order traditional upwind SBP operators have $b=4$. The numerical experiments in Section~\ref{sec:Computations_1D} indicate that the boundary-optimized operators can yield higher-than-expected convergence rates, in particular for odd orders.
\end{remark}

We denote by $p$ the interior order of accuracy. The construction of even-order diagonal-norm upwind (or traditional) SBP operators (up to eighth order) requires at least $p$ boundary points with boundary closures of accuracy $p/2$\revOne{~\cite{Strand94,MattssonAlmquistCarpenter13}}. Analogously, odd-order upwind (or traditional) SBP operators (up to ninth order) require at least $p-1$ boundary points with boundary closures of accuracy $(p-1)/2$.

\subsubsection{Construction procedure}\label{sec:procedure}

The construction starts by enforcing the required symmetry and SBP compatibility conditions while introducing sufficient degrees of freedom in the boundary closures to meet the accuracy conditions. The operators are derived using Maple (symbolic computations), although any comparable symbolic algebra software can be used.

The interior stencil for the $p^{\text{th}}$-order upwind SBP operator $D_{+}$ (and the corresponding $(p+1)^{\text{th}}$-order operator) is obtained from a central narrow stencil of order $p+2$ by subtracting
\[
\frac{h^{p+1}}{\alpha_p}\,(\Delta_{+}\Delta_{-})^{p/2+1}
\]
to remove the outermost grid point on the left. Here $\Delta_{\pm}$ denote the standard forward and backward difference operators, and $\alpha_p$ depends on the order. This yields a $(p+1)$-order accurate skewed stencil, which is the interior stencil for the $(p+1)^{\text{th}}$-order upwind SBP operator $D_{+}$. \revOne{For example, the construction of the interior stencil of the $3^{\text{rd}}$ order $D_{+}$ operator starts with the $4^{\text{th}}$ order narrow central stencil $\frac{1}{12h}[1,-8,0,8,-1]$ and subtracts $\frac{h^{3}}{12}\,(\Delta_{+}\Delta_{-})^{2} = \frac{1}{12h}[1, -4, 6, -4, 1]$ yielding the interior stencil $\frac{1}{12h}[0, -4, -6, 12, -2]$, a skewed stencil of order 3.} To obtain the interior stencil of the $p^{\text{th}}$-order operator, an additional subtraction is applied to remove the next outermost term on the left, yielding an off-centering by two grid points.

The matrix $Q_+$ is parameterized by the ansatz
\[
\left[ \begin {array}{cccccccccccc} 
q_{1,\,1}&q_{1,\,2}&\ldots&q_{1,\,\tilde{s}}& &&&&&&&\\
q_{2,\,1}   &\ddots    &     &q_{2,\,\tilde{s}}& &&&&&&&\\
\vdots      &              &      &\vdots   &  q_{\frac{s}{2}+1} & &&&& && \\
\vdots      &              &      &\vdots   & \vdots                   & &&&&   &&\\
q_{\tilde{s},\,1}  &q_{\tilde{s},\,2}&\ldots&q_{\tilde{s},\,\tilde{s}}& q_1    &\ddots & \ddots & &&&&\\         
 & &q_{n-\frac{s}{2}}&\ldots& q_{0} &q_{1}   &   \ldots&q_{\frac{s}{2}+1} &&&&\\      
             & & &\ddots &  &\ddots &&&\ddots &&&\\   
              && & &         q_{n-\frac{s}{2}}&\ldots& q_{0}&q_{1}&\ldots  &   q_{\frac{s}{2}+1}&&\\
           &&& & &            \ddots&                   & \ddots              & q_{\tilde{s},\,\tilde{s}}      &q_{\tilde{s}-1,\,s}& \ldots &q_{1,\,\tilde{s}}\\
                        &&&& &  &            \ddots            &               &q_{\tilde{s},\,\tilde{s}-1}    &\ddots      &            &q_{1,\,\tilde{s}-1}\\
                        &&&& & &                            &   q_{n-\frac{s}{2}}            &\vdots           &               &            &\vdots \\
                        &&&& & &                             &              &\vdots           &               &            &\vdots \\
                         &&&& & &                     &                     &q_{\tilde{s},\,1}   &\ldots &        &q_{1,\,1}\\
\end {array} \right].
\]
\revOne{In the ansatz above, $s$ is the width of the repeating interior stencil, $n$ is the offset, and $\tilde{s}\times\tilde{s}$ is the size of the boundary closure block. In the repeating interior stencil, the unknown $q_0$ is the diagonal entry of the matrix.}
In this work, boundary-optimized upwind SBP operators are derived for $p=2,4,6,8$ (even orders) and $p=3,5,7,9$ (odd orders). The first two grid points adjacent to each boundary are non-equidistant. For even orders, $s=p$ and $n=1$, while for odd orders $s=p-1$ and $n=2$ in the above structure. To obtain $s/2$-accurate boundary closures in both $D_+$ and $D_-$, at least $s$ boundary points are required, i.e., $\tilde{s}\ge s$. The operators derived here use the minimal choice $\tilde{s}=s$. The norm matrix is parameterized as
\[
H = h\,\mathrm{diag}(h_1,\ldots,h_s,1,\ldots,1,h_s,\ldots,h_1).
\]
The grid parameters $d_1$ and $d_2$ in \eqref{eq:gridpoints} are included among the unknowns.

The unknown parameters in ${\bf x}$, $Q_+$, and $H$ are determined by enforcing the accuracy requirements in Definition~\ref{definition:accuracy}, together with the SBP symmetry constraints. It is essential that $H$ remains positive definite. For $p>3$, free parameters remain after enforcing the accuracy and symmetry conditions. The fourth- and fifth-order operators admit four free parameters, the sixth- and seventh-order operators admit six, and the eighth- and ninth-order operators admit eleven.

When selecting these free parameters, two additional objectives are enforced:
\begin{enumerate}
\renewcommand\labelenumi{\theenumi)}
\item The dissipation matrix $S=\tfrac12(Q_+ + Q_+^T)$ is required to be negative semi-definite.
\item Boundary accuracy is optimized by minimizing the leading-order error norm in Definition~\ref{definition:error}.
\end{enumerate}
\revOne{The resulting constrained non-convex, nonlinear optimization problems are solved using the \texttt{Minimize} routine in Maple. The default convergence settings of \texttt{Minimize} are used; no user-specified optimization tolerance is imposed. All parameter optimizations are performed on the reference grid with $h=1$ and $m=21$ grid points described above. Multiple local minima may exist, and global optimality cannot be guaranteed. Order-dependent box constraints are imposed on the free grid and operator parameters to restrict the numerical search to reasonable boundary closures.}

\revOne{For example, considering the eighth-order operator we use
\[
d_1\in[0.1,0.7],\qquad d_2\in[0.3,1.2],
\]
while the remaining free $q$-parameters are restricted to $[-2,2]$.
The objective combines the two leading error vectors. 
For the eighth-order operator the two relevant error vectors \eqref{eq:error_vec} are ${\bf e}_{\pm,(q)}$ with $q=5,6$ obtained with $\hat Q=Q_\pm$.
Minimizing the cost function
\[
J_8=
\|{\bf e}_{+,(5)}+{\bf e}_{+,(6)}\|_h^2+
\|{\bf e}_{-,(5)}+{\bf e}_{-,(6)}\|_h^2.
\]
results in the boundary-optimized eighth-order operator.}

\revOne{Since $h=1$ during the optimization, this is exactly the sum of the corresponding Euclidean inner products used in the Maple implementation. This combined-error objective follows the optimization strategy used in \cite{MattssonAlmquistCarpenter13}. Other objective functions, such as the sum of the individual error norms used in later boundary-optimization work \cite{MattssonAlmquistWeide2018}, are also possible. The spectral radius is not included explicitly in the present objective function. After the free parameters have been fixed, the resulting boundary-optimized upwind SBP operators contain no additional tunable parameters.}

\section{Hyperbolic 1D system}\label{sec:1-D_1}

Next, we demonstrate how the upwind SBP operators, in combination with the projection method, are used to construct stable semi-discrete approximations of one-dimensional hyperbolic systems.

\subsection{Continuous problem}\label{sec:1-D_11}

A stability estimate is derived using an energy method. Consider the hyperbolic system
\begin{equation}\label{eq:PDE2_1D} 
\begin{array}{lll}
{\bf C}{\bf u}_{t}+{\bf A}{\bf u}_x=\revOne{\mathbf{0}},   & x \in[-1,\,1], &  t\ge 0,\\[0.1cm]
u^{(1)}(-1,t) = 0 , & &  t\ge 0,\\[0.1cm]
u^{(1)}(1,t)=0, & &  t\ge 0,\\[0.1cm]
{\bf u}(x,0)={\bf f}(x),& x \in[-1,\,1], &  \\
 \end{array}
\end{equation}
where ${\bf C}$ is diagonal and positive definite and ${\bf f}(x)$ denotes the initial data. The unknown and flux matrix are
\begin{equation}\label{eq:Matrix_problem1}
{\bf u}=\begin{bmatrix}
u^{(1)}\\
u^{(2)}\\
\end{bmatrix},\qquad
{\bf A}=\begin{bmatrix}
\alpha&1 \\
1&0\\
\end{bmatrix},
\qquad \alpha \in \mathbb{R}.
\end{equation}

Multiplying the PDE in \eqref{eq:PDE2_1D} by ${\bf u}^T$, integrating by parts, and adding the transpose yields the energy identity
\begin{equation}\label{eq:energy_continuous}
\frac{d}{dt}\|{\bf u}\|_C^2
=
-{\bf u}_r^T{\bf A}{\bf u}_r + {\bf u}_l^T{\bf A}{\bf u}_l
=
-u^{(1)}_{r}(\alpha u^{(1)}_{r}+ 2u^{(2)}_{r})  + u^{(1)}_{l}(\alpha u^{(1)}_{l} + 2u^{(2)}_{l}),
\end{equation}
where ${\bf u}_{l,r}^T=[u^{(1)}_{l,r}\ \ u^{(2)}_{l,r}]$ denote the boundary traces at $x=-1$ and $x=1$, respectively. Imposing $u^{(1)}(\pm 1,t)=0$ gives
\[
\frac{d}{dt}\|{\bf u}\|_C^2 = 0,
\]
and the problem is therefore energy stable. In the following, we assume $\alpha \ge 0$.

\subsection{Semi-discrete problem}\label{sec:1-D_12}

\revOne{To distinguish the continuous and discrete formulations, boldface is used for the continuous vector and matrix quantities in Section~\ref{sec:1-D_11}, while their discrete counterparts in this section are written in plain type.}

\revThree{To construct the SBP-projection discretization, we first introduce the unprojected SBP spatial operator together with the discrete boundary constraints,}
\begin{equation}\label{eq:SBP3} 
\begin{array}{ll}
Cv_t + D_x v = 0, & t \ge 0,\\[0.1cm]
(e^{(1)} \otimes e_1^T)  v = 0, & t \ge 0,\\[0.1cm]
(e^{(1)} \otimes e_m^T) v = 0, & t \ge 0,\\[0.1cm]
v = f, & t = 0.
\end{array}
\end{equation}
\revThree{Equation~\eqref{eq:SBP3} is not solved as an overdetermined semi-discrete system; rather, it specifies the spatial operator and the boundary constraint used to define the projection.}
\revThree{where ${v} = [v^{(1)}\ , \ v^{(2)}] = [v_1^{(1)}\ldots v_m^{(1)}\ , \ v_1^{(2)}\ldots v_m^{(2)}]$, $f = [f^{(1)}\ , \ f^{(2)}] = [f^{(1)}(x_1)\ldots f^{(1)}(x_m)\ , \ f^{(2)}(x_1)\ldots f^{(2)}(x_m)]$,} and $e^{(1)} = [1 \ , \ 0]$. The boundary operator is defined as
\[
L =
\begin{bmatrix}
e^{(1)} \otimes e_1^T \\
e^{(1)} \otimes e_m^T
\end{bmatrix},
\]
and the discrete spatial operator is given by
\[
D_x =
\begin{bmatrix}
\alpha D_- & D_+ \\
D_- & 0
\end{bmatrix}.
\]

\revThree{The SBP-projection approximation of the continuous problem \eqref{eq:PDE2_1D} is then}
\begin{equation}\label{eq:SBP_Projection2} 
\begin{array}{ll}
v_t + P C^{-1} D_x P v = 0, & t\ge 0,\\
v=f, & t=0,
\end{array}
\end{equation}
where the projection operator is given by
\[
P = I - \bar{H}^{-1}L^T\left(L\bar{H}^{-1}L^T\right)^{-1}L,
\qquad \bar{H}=H_2 C.
\]
The SBP structure follows from
\[
H_2 D_x=
\begin{bmatrix}
\alpha Q_-+\alpha\frac{B}{2} & Q_+ + \frac{B}{2}\\
Q_-+\frac{B}{2} & 0
\end{bmatrix},
\]
where \revOne{$H_2 = I_2\otimes H$ with $I_2$ as the $2\times 2$ identity matrix and} $B$ is defined in \eqref{eq:Def_matrix}. Hence,
\begin{equation}\label{eq:BBar}
H_2D_x+(H_2D_x)^T=\bar{B}=
\begin{bmatrix}
-2\alpha S+\alpha B & B\\
B & 0
\end{bmatrix}.
\end{equation}

\begin{remark}
For central-difference SBP operators, replacing $D_{\pm}$ by $D_1$ in $D_x$ yields
\[
H_2D_x+(H_2D_x)^T=\bar{B}=
\begin{bmatrix}
\alpha B & B\\
B & 0
\end{bmatrix}.
\]
\end{remark}

The identity $\bar{H}P=P^T\bar{H}$ is used in the subsequent discrete energy analysis. Introducing the projected variables $w=Pv$, multiplying \eqref{eq:SBP_Projection2} by $v^T\bar{H}$, and adding the transpose gives
\begin{equation}\label{eq:Stability2_P} 
\begin{aligned}
\frac{d}{dt}\| v \|_{\bar{H}}^2
&= - (Pv)^T\left(H_2D_x+(H_2D_x)^T\right)(Pv)\\
&= - w^T \bar{B} w\\
&= -w^{(1)}_m\left(\alpha w^{(1)}_m+2w^{(2)}_m\right)
+w^{(1)}_1\left(\alpha w^{(1)}_1+2w^{(2)}_1\right)
+2\alpha (w^{(1)})^T S w^{(1)}.
\end{aligned}
\end{equation}
Since $Lw=0$, the projected solution satisfies the boundary conditions exactly, i.e., $w^{(1)}_1=w^{(1)}_m=0$, and therefore
\[
\frac{d}{dt}\| v \|_{\bar{H}}^2 = 2\alpha (w^{(1)})^T S w^{(1)} \le 0,
\]
where the inequality follows from the defining property $S=S^T\le 0$ of the upwind SBP operators. For central-difference SBP operators, $S=0$, i.e., no dissipation is introduced for $\alpha>0$.
\revOne{
\begin{remark}
For the 1D tests considered herein, we use homogeneous boundary conditions. The projection formulation can also be extended to inhomogeneous boundary data by introducing a lifting $\tilde g(t)$ satisfying $L\tilde g(t)=g(t)$ and projecting the homogeneous remainder. The corresponding construction and stability analysis are given in \cite{MattssonOlsson2018}; we do not pursue this extension here.
\end{remark}
}

\section{Computations in 1D}\label{sec:Computations_1D}

We investigate convergence using a smooth analytical solution for $\alpha=0$ and compare the proposed boundary-optimized upwind SBP operators with the traditional upwind SBP operators presented in \cite{Mattsson17}. For completeness, we also compare with traditional central-difference SBP operators $D_1$ and the corresponding boundary-optimized $D_1$ operators derived in \cite{MattssonAlmquistWeide2018}, for orders up to twelve.

   \revOne{The traditional central-difference SBP operators are included primarily as reference discretizations, particularly for the the one-dimensional non-smooth test, to provide a comparison with the proposed upwind operators. We emphasize that such nondissipative central schemes are not intended to represent a preferred discretization for nonlinear hyperbolic problems. Instead, the traditional central and upwind schemes, without artificial dissipation, provide a deliberately stringent comparison that allow us to assess the accuracy and stability properties of the proposed boundary-optimized upwind operators independent of any additional dissipation mechanism.}

Robustness is further examined by introducing convection through $\alpha=3$, which leads to the development of non-smooth solution features after boundary interaction. In all numerical experiments, we set $C$ equal to the identity matrix.

\subsection{Spectral radius}\label{sec:Interior}

We first compare the spectral radius of the SBP-projection approximation \eqref{eq:SBP_Projection2} for $\alpha=0$. The discretization matrix is given by
\[
M=-PC^{-1}D_xP.
\]
Table~\ref{table:Spectral} reports the spectral radius of the normalized matrix $hM$, comparing the traditional upwind SBP operators from \cite{Mattsson17} with the proposed boundary-optimized upwind SBP operators. The eighth-order boundary-optimized upwind operator exhibits a comparatively large spectral radius and is therefore expected to be less efficient for explicit time integration. \revOne{The spectral radius was not included explicitly in the optimization objective. The eighth-order constrained optimization proved considerably more difficult than the lower-order cases. Several feasible parameter choices were investigated, but we did not obtain a closure that retained comparably small boundary errors while substantially reducing the spectral radius. The reported eighth-order closure should therefore not be interpreted as being jointly optimal with respect to both boundary accuracy and explicit time stepping efficiency. A systematic multi-objective optimization that includes the spectral radius is left for future work.}

\begin{table}
\caption{The spectral radius of $h\,M$ comparing traditional upwind and the boundary-optimized upwind SBP operators.} \label{table:Spectral}
\label{tab:spectra}
\centering
\begin{tabular}{l c c c}
\toprule
 Operator & $p=4$  & $p=6$ & $p=8$  \\
\midrule
Traditional upwind, order $p$ & 2.6652 &3.1060& 2.3706 \\
Traditional upwind, order $p+1$     & 1.6234  & 1.7452&    2.0260 \\
Optimimized upwind, order $p$ & 2.6310 &2.5732& 17.6784\\
Optimimized upwind, order $p+1$ & 1.9280 &2.4187& 2.7045\\
 \bottomrule
\end{tabular}
\end{table}

\subsection{Convergence}\label{sec:Convergence}

The accuracy of the SBP-projection approximation \eqref{eq:SBP_Projection2} is verified against an analytical solution in which a narrow Gaussian pulse propagates and reflects at the boundaries. In this section we set $\alpha=0$, corresponding to pure wave propagation without convection.

Let $u(t_*)$ and $u_{(m)}$ denote the analytical and numerical solutions at time $t=t_*$, evaluated on a grid with $m$ points. The error vector is defined as
\[
err_{(m)} = u(t_*) - u_{(m)},
\]
and the discrete $\ell^2$-norm of the error is
\[
\|err_{(m)}\|_h=\sqrt{h}\,\|err_{(m)}\|.
\]
\revOne{This scaled Euclidean norm is used consistently for the one-dimensional convergence tests. Since the non-equidistant region contains only a fixed number of boundary points and the diagonal SBP weights remain uniformly bounded, this norm and the corresponding SBP norm are equivalent and yield the same asymptotic convergence rate.}
The observed convergence rate, denoted by $q$, is computed as
\begin{equation}\label{eq:noggr}
q=\log_{10}\left(\frac{\|err_{(m)}\|_h}{\|err_{(n)}\|_h}\right)/\log_{10}\left(\frac{n}{m}\right)\;,
\end{equation}
where $m$ and $n$ denote two successive grid resolutions.

The Gaussian profiles
\[
\theta^{(1)}(x,t)=\exp{\left(-\left(\frac{x-t}{r_{*}}\right)^2\right)}\quad,\quad
\theta^{(2)}(x,t)=-\exp{\left(-\left(\frac{x+t}{r_{*}}\right)^2\right)}\;,
\]
are introduced, where $r_{*}$ defines the pulse width. The initial data are set to
\begin{equation}\label{eq:initial}
u^{(1)}(x,0)=\theta^{(1)}(x,0)-\theta^{(2)}(x,0)\,,\quad u^{(2)}(x,0)=\theta^{(1)}(x,0)+\theta^{(2)}(x,0)\;.
\end{equation}
On a domain of length $L=2$, the analytical solution to \eqref{eq:PDE2_1D} with $\alpha=0$ and initial data \eqref{eq:initial} is given by
\[
u^{(1)}(x,t)=\theta^{(2)}(x,L-t)-\theta^{(1)}(x,L-t)\,,\quad u^{(2)}(x,t)=\theta^{(1)}(x,L-t)+\theta^{(2)}(x,L-t)\;,
\]
after the Gaussian pulses have been reflected at the boundaries, i.e., for $t \in  [L-L/8,\, L+L/8]$.

We set $r_{*}=0.1$ and integrate to $t=t_*=1.8$ using the classical four-stage, fourth-order Runge--Kutta method~\cite{butcher2008numerical} with CFL number $0.05$, i.e., $k =0.05\,h$, so that temporal errors are negligible compared with spatial errors. Table~\ref{table:convergence_standard} reports convergence for traditional central-difference SBP operators of orders $6$--$12$, and Table~\ref{table:convergence_standard_optimal} reports the corresponding results for boundary-optimized central-difference SBP operators. Tables~\ref{table:convergence_even_optimal} and \ref{table:convergence_odd_optimal} report convergence for the boundary-optimized upwind SBP operators, while Tables~\ref{table:convergence_even_traditional} and \ref{table:convergence_odd_traditional} report results for the traditional upwind SBP operators.

\begin{table}[H]
\centering
\caption{$\log(\ell^2\ \text{errors})$ and convergence rates for \eqref{eq:SBP_Projection2} using even-order ($6^\text{th}$--$12^\text{th}$) traditional central-difference SBP operators from \cite{MattssonNordstrom04}, CFL=0.05.}
\label{table:convergence_standard}
\begin{tabular}{lllllllll}
\toprule
  $m$  & $\log{l^2}_{(6)}$ & $q$ & $\log{l^2}_{(8)}$& $q$& $\log{l^2}_{(10)}$& $q$ & $\log{l^2}_{(12)}$& $q$ \\
\midrule
$ 51$ & -1.05 & --- & -0.83 & --- & -0.72 & --- & -0.58 & --- \\
$101$ & -1.74 & 2.28 & -2.09 & 4.18 & -1.98 & 4.19 & -1.53 & 3.16 \\
$201$ & -2.87 & 3.77 & -3.47 & 4.60 & -3.43 & 4.82 & -3.22 & 5.62 \\
$401$ & -4.27 & 4.66 & -4.83 & 4.49 & -5.24 & 6.01 & -5.46 & 7.43 \\
 \bottomrule
\end{tabular}
\end{table}

\begin{table}[H]
\centering
\caption{$\log(\ell^2\ \text{errors})$ and convergence rates for \eqref{eq:SBP_Projection2} using even-order ($6^\text{th}$--$12^\text{th}$) boundary-optimized central-difference SBP operators from \cite{MattssonAlmquistWeide2018}, CFL=0.05.}
\label{table:convergence_standard_optimal}
\begin{tabular}{lllllllll}
\toprule
  $m$  & $\log{l^2}_{(6)}$ & $q$ & $\log{l^2}_{(8)}$& $q$& $\log{l^2}_{(10)}$& $q$ & $\log{l^2}_{(12)}$& $q$ \\
\midrule
$ 51$ & -1.18 & --- & -1.55 & --- & -1.85 & ---  & -1.97 & --- \\
$101$ & -2.70 & 5.05 & -3.59 & 6.77 & -4.27 & 8.04 & -4.42 & 8.16 \\
$201$ & -4.47 & 5.88 & -5.94 & 7.80 & -7.20 & 9.72 & -7.55 & 10.40 \\
$401$ & -6.25 & 5.90 & -8.28 & 7.77 & -9.51 & 7.68 & -10.73 & 10.57 \\
 \bottomrule
\end{tabular}
\end{table}

\begin{table}[H]
\centering
\caption{$\log(\ell^2\ \text{errors})$ and convergence rates for \eqref{eq:SBP_Projection2} using even-order ($4^\text{th}$--$8^\text{th}$) boundary-optimized upwind SBP operators derived herein, CFL=0.05.}
\label{table:convergence_even_optimal}
\begin{tabular}{lllllll}
\toprule
  $m$  & $\log{l^2}_{(4)}$ & $q$ & $\log{l^2}_{(6)}$& $q$& $\log{l^2}_{(8)}$& $q$ \\
\midrule
$ 51$ & -0.39 & --- & -1.19 & --- & -2.10 & --- \\
$101$ & -1.17 & 2.60 & -2.78 & 5.30 & -4.27 & 7.20 \\
$201$ & -2.32 & 3.82 & -4.56 & 5.93 & -6.65 & 7.88 \\
$401$ & -3.52 & 3.98 & -6.37 & 5.99 & -9.05 & 7.98 \\
 \bottomrule
\end{tabular}
\end{table}

\begin{table}[H]
\centering
\caption{$\log(\ell^2\ \text{errors})$ and convergence rates for \eqref{eq:SBP_Projection2} using even-order ($4^\text{th}$--$8^\text{th}$) traditional upwind SBP operators from \cite{Mattsson17}, CFL=0.05.}
\label{table:convergence_even_traditional}
\begin{tabular}{lllllll}
\toprule
  $m$  & $\log{l^2}_{(4)}$ & $q$ & $\log{l^2}_{(6)}$& $q$& $\log{l^2}_{(8)}$& $q$ \\
\midrule
$ 51$ & -0.42 & --- & -0.65 & --- & -1.28 & --- \\
$101$ & -1.20 & 2.60 & -2.70 & 6.78 & -3.64 & 7.86 \\
$201$ & -2.34 & 3.78 & -5.21 & 8.36 & -6.41 & 9.20 \\
$401$ & -3.53 & 3.96 & -6.80 & 5.27 & -9.40 & 9.92 \\
 \bottomrule
\end{tabular}
\end{table}

\begin{table}[H]
\centering
\caption{$\log(\ell^2\ \text{errors})$ and convergence rates for \eqref{eq:SBP_Projection2} using odd-order ($5^\text{th}$--$9^\text{th}$) boundary-optimized upwind SBP operators derived herein, CFL=0.05.}
\label{table:convergence_odd_optimal}
\begin{tabular}{lllllll}
\toprule
  $m$  & $\log{l^2}_{(5)}$ & $q$ & $\log{l^2}_{(7)}$& $q$& $\log{l^2}_{(9)}$& $q$ \\
\midrule
$ 51$ & -1.22 & --- & -2.01 & --- & -2.48 & --- \\
$101$ & -2.88 & 5.50 & -4.31 & 7.67 & -5.70 & 10.72 \\
$201$ & -4.69 & 6.00 & -6.71 & 7.96 & -8.82 & 10.36 \\
$401$ & -6.50 & 6.03 & -9.11 & 7.97 & -11.76 & 9.77 \\
 \bottomrule
\end{tabular}
\end{table}

\begin{table}[H]
\centering
\caption{$\log(\ell^2\ \text{errors})$ and convergence rates for \eqref{eq:SBP_Projection2} using odd-order ($5^\text{th}$--$9^\text{th}$) traditional upwind SBP operators from \cite{Mattsson17}, CFL=0.05.}
\label{table:convergence_odd_traditional}
\begin{tabular}{lllllll}
\toprule
  $m$  & $\log{l^2}_{(5)}$ & $q$ & $\log{l^2}_{(7)}$& $q$& $\log{l^2}_{(9)}$& $q$ \\
\midrule
$ 51$ & -0.78 & --- & -1.12 & --- & -1.36 & --- \\
$101$ & -2.27 & 4.94 & -2.79 & 5.54 & -3.53 & 7.20 \\
$201$ & -3.89 & 5.37 & -5.15 & 7.84 & -6.20 & 8.88 \\
$401$ & -5.43 & 5.14 & -7.28 & 7.09 & -8.90 & 8.98 \\
 \bottomrule
\end{tabular}
\end{table}

The expected convergence rates are $3,4,5,6,$ and $7$ for even-order SBP operators of orders $4,6,8,10,$ and $12$, respectively. For odd-order SBP operators of orders $5,7,$ and $9$, the expected convergence rates are $3,4,$ and $5$. At coarse and moderate resolutions, the boundary-optimized operators---in particular the proposed upwind SBP operators---exhibit substantially higher observed convergence rates. The boundary-optimized odd-order upwind SBP operators demonstrate especially strong accuracy at marginal resolution.

\revTwo{The results from Tables~\ref{table:convergence_standard} to \ref{table:convergence_odd_traditional} show that the boundary-optimized upwind SBP operators generally provide substantially higher observed convergence rates. 
For example, the observed convergence rates of the fifth-, seventh-, and ninth-order upwind operators for the boundary-optimized upwind SBP operators reach approximately 6.0, 8.0, and 10.4, respectively, compared to the approximately 5.4, 7.8, and 8.9 for the corresponding traditional central (or upwind) SBP operators. The improvement is particularly pronounced for the ninth-order operator, for which the optimized operator maintains convergence rates close to 10 over several successive refinements. The error magnitudes show a similar trend. The largest improvements occur at the lower and intermediate resolutions, consistent with the objective of reducing the leading boundary truncation-error contributions. At the finest resolutions, the differences become less pronounced as higher-order truncation-error terms become increasingly important.}

\revTwo{This surprising reduction in the error at moderate resolution may arise from suppressing the dominant boundary truncation error terms due to the boundary-optimization procedure. As such, once those terms become sufficiently small, higher-order terms and/or interior errors become relevant, producing an apparent superconvergent regime over the resolutions considered. However, more analysis and investigation is necessary to quantify the improved convergence rates of these boundary-optimized SBP operators.}

\subsection{Non-smooth boundary interaction}\label{sec:Interaction}

We next set $\alpha=3$, which yields predominantly right-going waves. 
After interaction with the right boundary, the solution develops non-smooth features. The initial data are again given by \eqref{eq:initial}.
\revOne{This test case configuration is used to assess the robustness and behavior of the discrete solution following a boundary interaction that develops non-smooth features, comparing the newly described boundary-optimized upwind SBP operators against those available in the literature. The novel operators exhibit significantly less spurious ringing compared to other available, high-order SBP operators.}

\begin{figure}[!h]
\centering

\begin{subfigure}[b]{0.49\textwidth}
  \centering
  \includegraphics[width=\linewidth]{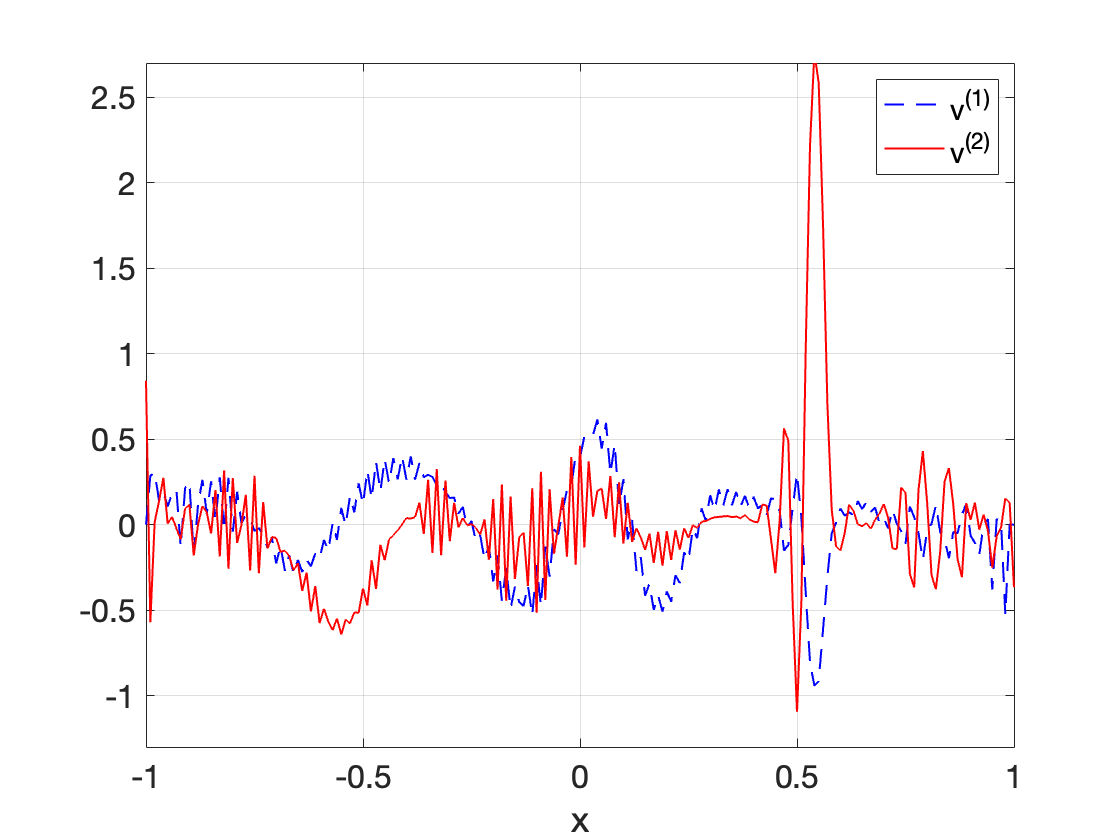}
  \caption{$m=201$, traditional SBP central $12^\text{th}$ order}
\end{subfigure}
\hfill
\begin{subfigure}[b]{0.49\textwidth}
  \centering
  \includegraphics[width=\linewidth]{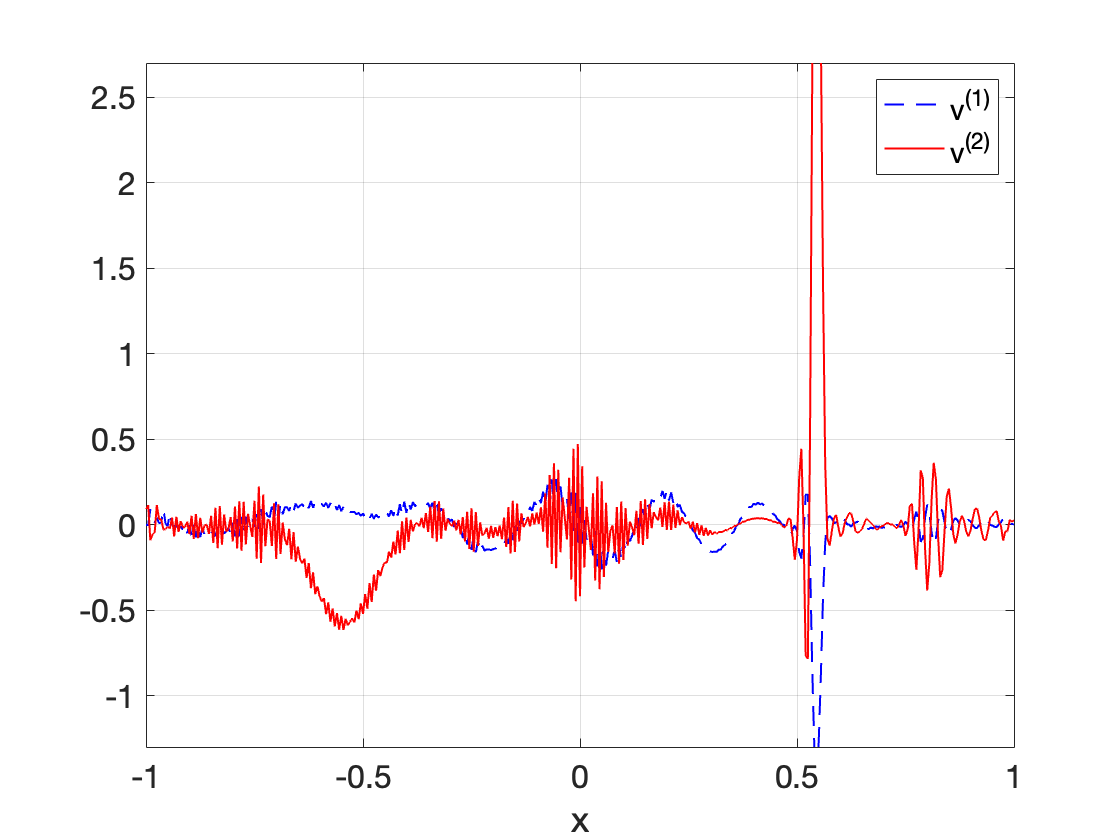}
  \caption{$m=401$, traditional SBP central $12^\text{th}$ order}
\end{subfigure}

\begin{subfigure}[b]{0.49\textwidth}
  \centering
  \includegraphics[width=\linewidth]{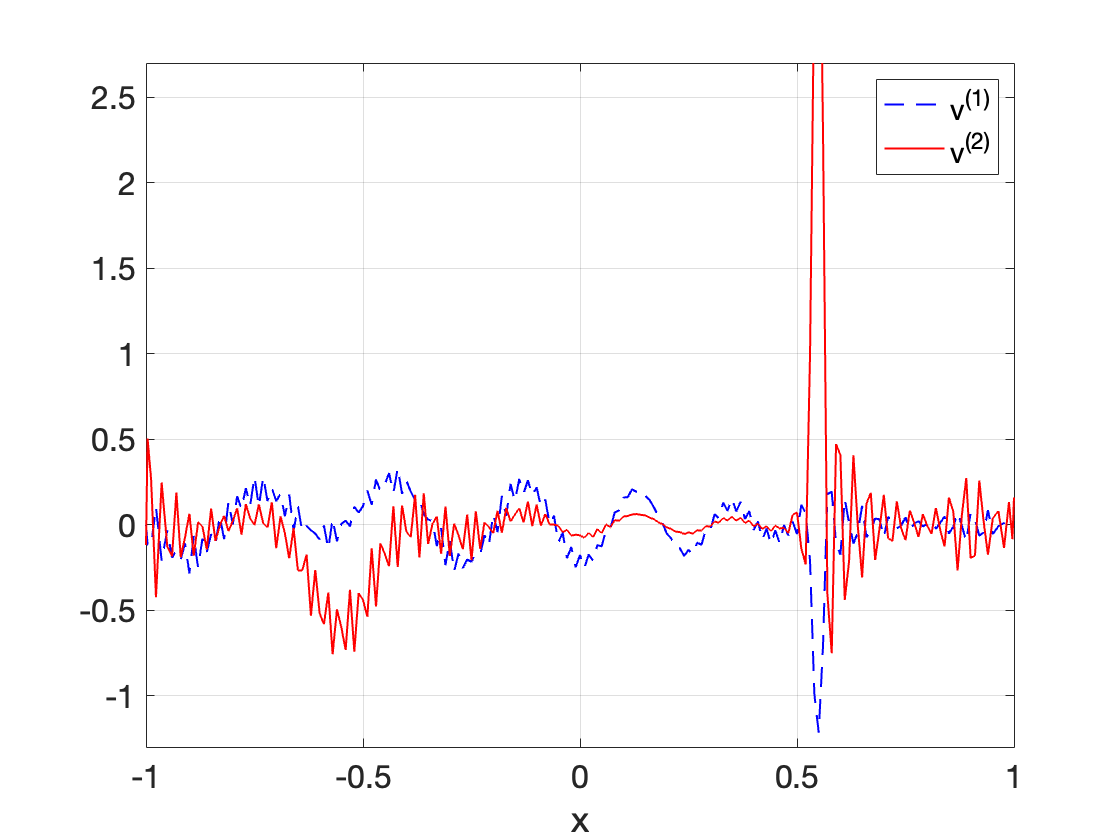}
  \caption{$m=201$, optimal SBP central $12^\text{th}$ order}
\end{subfigure}
\hfill
\begin{subfigure}[b]{0.49\textwidth}
  \centering
  \includegraphics[width=\linewidth]{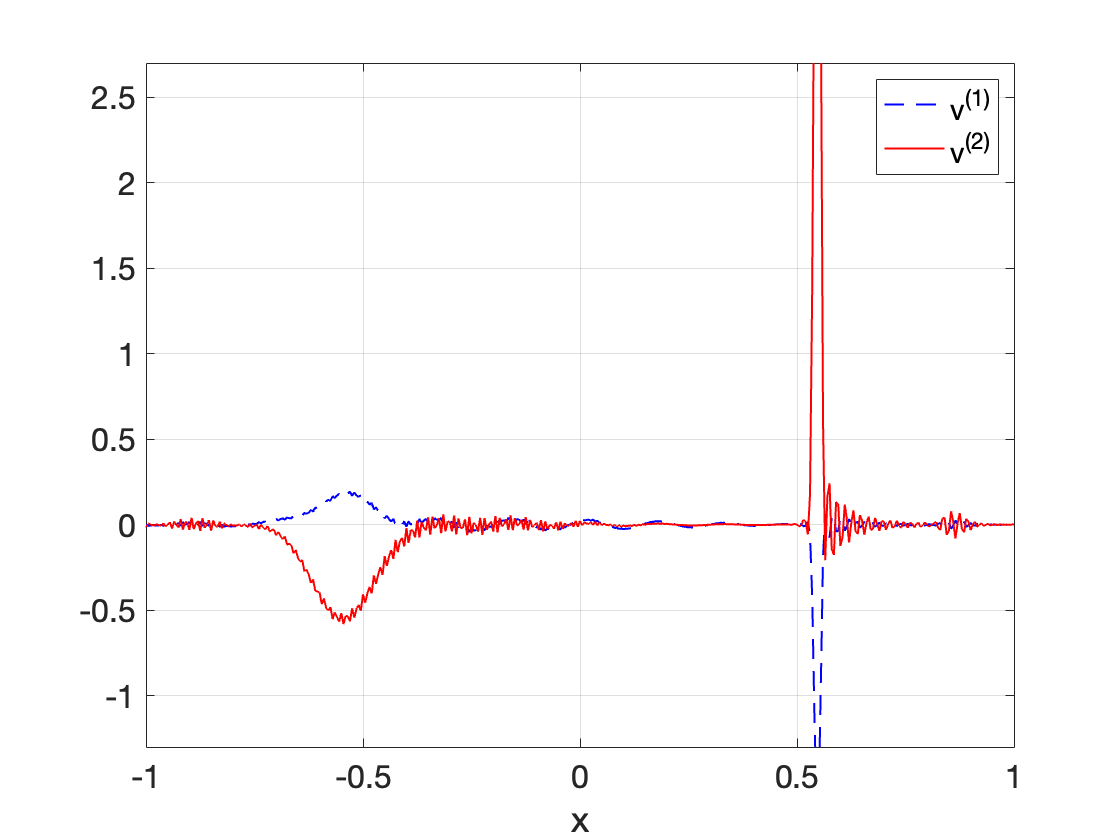}
  \caption{$m=401$, optimal SBP central $12^\text{th}$ order}
\end{subfigure}

\begin{subfigure}[b]{0.49\textwidth}
  \centering
  \includegraphics[width=\linewidth]{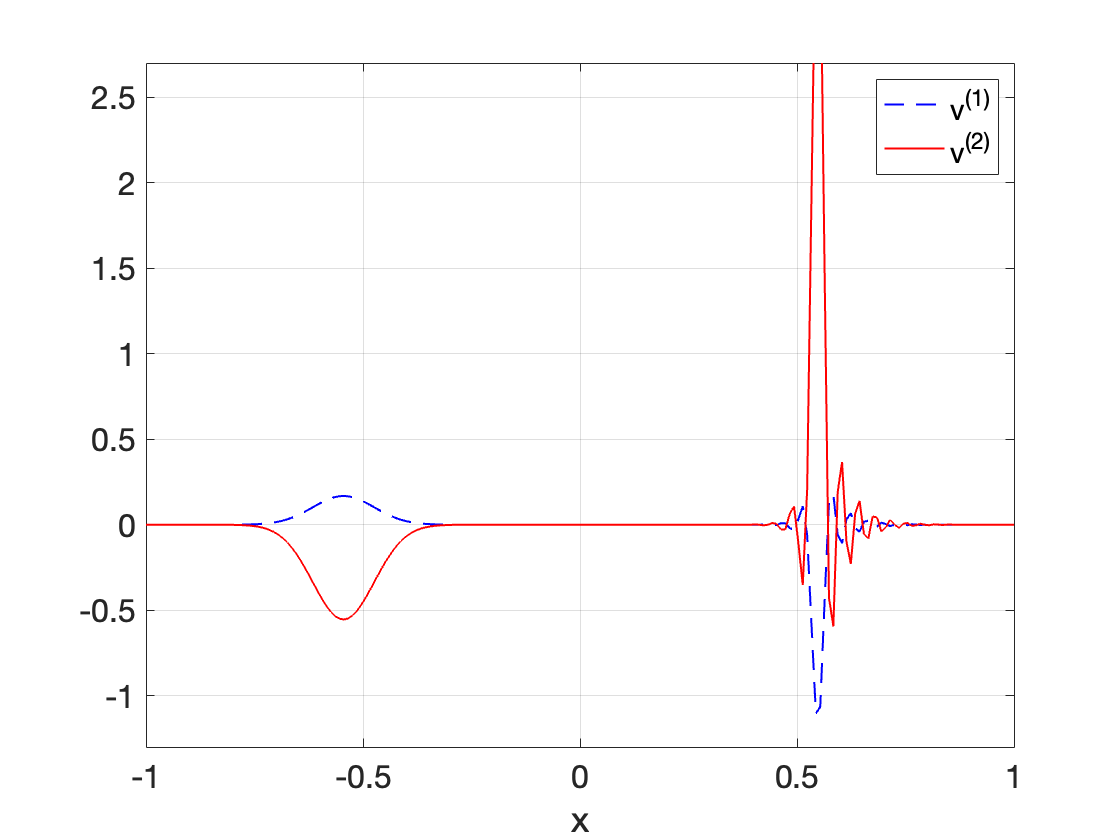}
  \caption{$m=201$, optimal SBP upwind $9^\text{th}$ order}
\end{subfigure}
\hfill
\begin{subfigure}[b]{0.49\textwidth}
  \centering
  \includegraphics[width=\linewidth]{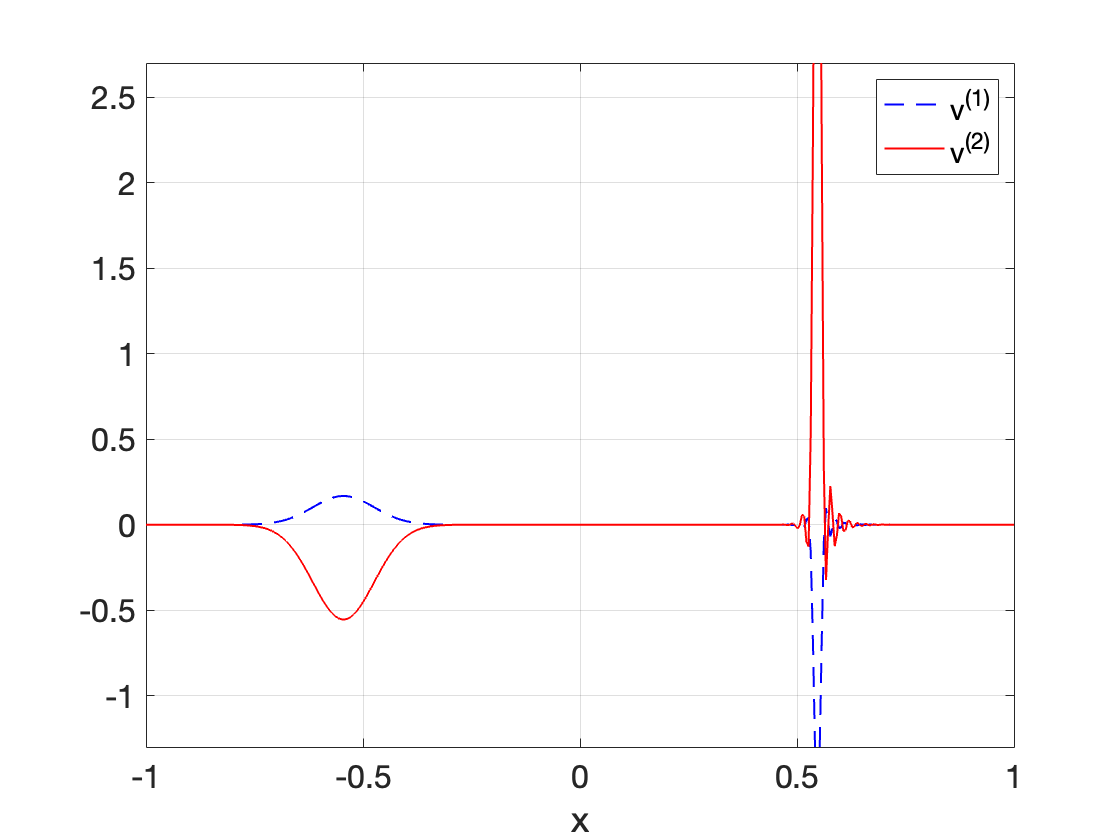}
  \caption{$m=401$, optimal SBP upwind $9^\text{th}$ order}
\end{subfigure}

\caption{Comparing central difference SBP (both traditional and boundary-optimized) against the boundary-optimized $9^\text{th}$-order upwind SBP operator for two grid resolutions. The simulations are run to $t=1.8$ with CFL=0.05 and $\alpha=3$.}
\label{fig:Convection}
\end{figure}

We set $r_{*}=0.1$ and integrate to $t=t_*=1.8$ using the classical four-stage, fourth-order Runge--Kutta method with CFL number $0.05$. Two grid resolutions are considered, $m=201$ and $m=401$. We compare the traditional $12^\text{th}$-order central SBP operator with the boundary-optimized $12^\text{th}$-order central SBP operator and the boundary-optimized $9^\text{th}$-order upwind SBP operator. The results are shown in Figure~\ref{fig:Convection}.

\section{Computations in 2D}\label{sec:Computations_2D}

Next, the newly derived upwind SBP operators are applied to nonlinear two-dimensional problems. To this end, the compressible Euler equations are considered,
\begin{equation}
  \label{eq:2d_comp_euler}
  \partial_t \begin{pmatrix} \rho \\ \rho v_1 \\ \rho v_2 \\ \rho e \end{pmatrix}
  + \partial_x \begin{pmatrix} \rho v_1 \\ \rho v_1^2 + p \\ \rho v_1 v_2 \\ (\rho e + p) v_1 \end{pmatrix}
  + \partial_y \begin{pmatrix} \rho v_2 \\ \rho v_1 v_2 \\ \rho v_2^2 + p \\ (\rho e + p) v_2 \end{pmatrix}
  = 0,
\end{equation}
describing the motion of an ideal gas with density $\rho$, velocity components $v_1$ and $v_2$, and total energy density $\rho e$. The pressure is given by
\begin{equation}
  p = (\gamma - 1) \left( \rho e - \frac{1}{2} \rho (v_1^2 + v_2^2) \right),
\end{equation}
where the ratio of specific heats is set to $\gamma = 1.4$. For all numerical experiments presented below, periodic boundary conditions are imposed on multi-block Cartesian meshes.

The spatial discretizations are implemented in Trixi.jl \cite{ranocha2022adaptive,schlottkelakemper2021purely}, where the boundary-optimized upwind SBP operators are provided through the package SummationByPartsOperators.jl \cite{ranocha2021sbp}. Time integration is performed using explicit Runge--Kutta methods from OrdinaryDiffEq.jl \cite{rackauckas2017differentialequations}; the specific schemes employed are described in the respective subsections. Visualization of the numerical solutions is carried out using ParaView \cite{ahrens2005paraview}.

\subsection{Isentropic vortex}

To examine the convergence properties of the boundary-optimized upwind SBP operators for a two-dimensional nonlinear problem, we consider the classical isentropic vortex test case introduced by Shu~\cite{shu1997essentially}. The analytical solution for given values of $x$, $y$, and $t$ is
\begin{equation}
\label{eq:isen_vort}
  \begin{aligned}
  \rho &= \left(1 - \frac{\varepsilon^2 (\gamma - 1) M^2}{8\pi^2}\exp\left(f(x,y,t)\right)\right)^{1 / (\gamma - 1)},
  \quad
  p = \frac{\rho^\gamma}{\gamma M^2},\\[0.1cm]
  v_1 &= 1 - \frac{\varepsilon (y-y_0)}{2\pi}\exp\left(\frac{f(x,y,t)}{2}\right),
  \quad
  v_2 = \frac{\varepsilon (x-x_0-t)}{2\pi}\exp\left(\frac{f(x,y,t)}{2}\right),
  \end{aligned}
\end{equation}
where
\[
f(x,y,t) = 1 - \left((x-x_0-t)^2 + (y-y_0)^2\right).
\]
Here, $M$ denotes the Mach number, 
$(x_0, y_0)$ specifies the initial vortex position, and $\varepsilon$ is the nondimensional circulation controlling the vortex strength. The initial condition is obtained by evaluating the vortex solution at $t=0$.

\begin{figure}[H]
\centering

\begin{subfigure}[b]{0.47\textwidth}
  \centering
  \includegraphics[width=\linewidth, trim=0 0 0 0, clip]{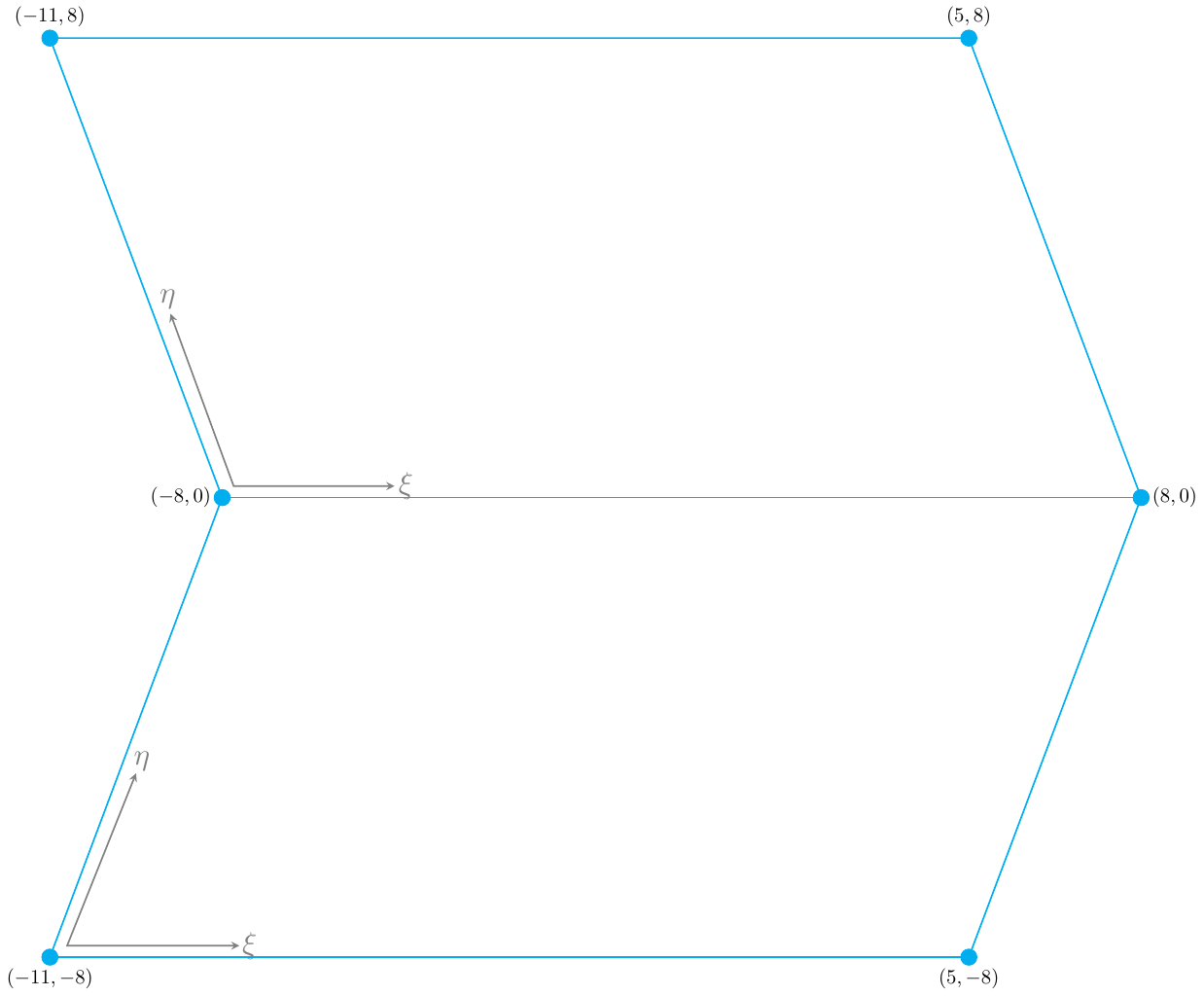}
  \caption{Chevron shaped domain with two blocks}
\end{subfigure}
\hspace*{0.5cm}
\begin{subfigure}[b]{0.47\textwidth}
  \centering
  \includegraphics[width=\linewidth, trim=0 0 0 0, clip]{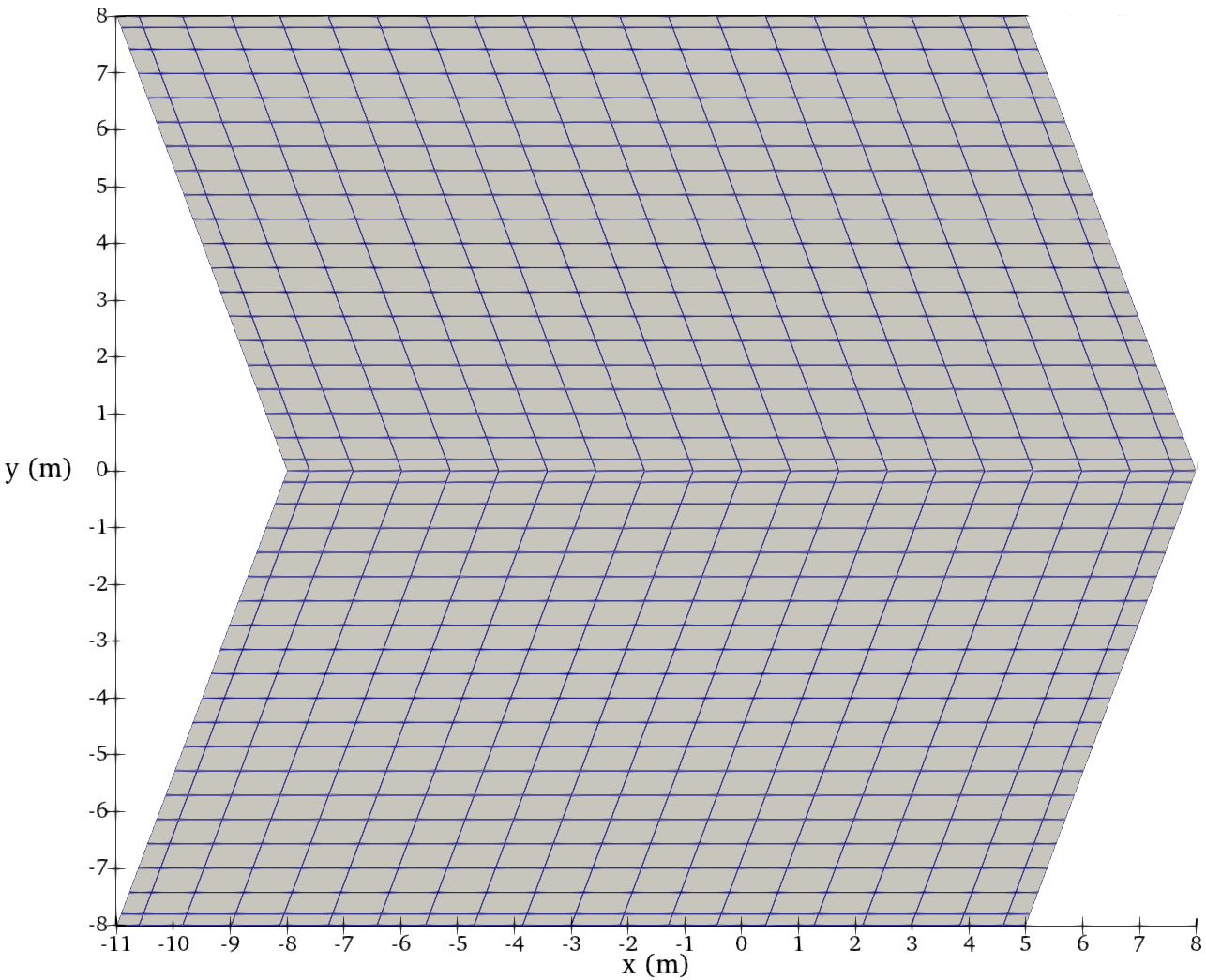}
  \caption{Discretization of domain with $21\times 21$ nodes in each block}
\end{subfigure}

\caption{Domain for the isentropic vortex convergence test and an example discretization of the two block domain with a non-uniform grid.}
\label{fig:domain}
\end{figure}

The isentropic vortex is simulated on a two-dimensional chevron-shaped domain depicted in Figure~\ref{fig:domain}. The domain is decomposed into two blocks; an example configuration using $21 \times 21$ grid points in each block is shown in Figure~\ref{fig:domain}. The node distribution associated with the boundary-optimized discretization is non-uniform. A mapping between computational and physical coordinates is applied in each block, introducing metric terms into the discretization; see, e.g., \cite{MattssonAlmquistCarpenter13}.

The initial vortex position is set to $(x_0, y_0) = (0, -0.5)$ such that the block interface, where discontinuities in the metric terms occur, intersects the vortex~\cite{MattssonAlmquistCarpenter13}. The vortex propagates from left to right across the slanted periodic boundaries. The isentropic vortex provides an analytical solution to the compressible Euler equations provided that \revOne{the domain is sufficiently large so that the discrete solution at the final simulation time does not interact with} the periodic boundaries\revTwo{~\cite{spiegel2015survey}}. Accordingly, for the convergence study, the final time is chosen as $t = 1$, ensuring that boundary effects do not influence the solution. For the convergence experiments, the Mach number is set to $M = 0.5$ and the vortex strength to $\varepsilon = 5$.

For the upwind discretizations, the local Lax--Friedrichs flux vector splitting for the compressible Euler equations is employed \cite{ranocha2025robustness}. 
\revTwo{The local Lax-Friedrichs splitting is, in general, consistent, conservative, and high-order accurate \cite{ranocha2025robustness, rydin2018high}; however, local linear and nonlinear stability of the local Lax-Friedrichs splitting has only been proven for the Burgers' equation \cite{ranocha2025robustness}. Nevertheless, the local Lax-Friedrichs splitting has been successfully applied to compressible Euler flows \cite{rydin2018high, duan2025active}.}
Coupling for standard SBP discretizations is performed using the local Lax--Friedrichs numerical flux function, which is equivalent to an upwind simultaneous approximation term (SAT).

Time integration of the semidiscrete system is carried out using the fourth-order accurate, nine-stage Runge--Kutta method proposed by Ranocha et al.~\cite{ranocha2021optimized}, combined with error-based step-size control and a sufficiently small tolerance of $10^{-12}$. This tolerance ensures that temporal errors are negligible and that the numerical approximation is dominated by spatial discretization errors.

Convergence results with respect to the density $\rho$ are presented for the boundary-optimized upwind SBP operators of odd and even orders in Tables~\ref{tab:convergence_odd_opt_upwind} and \ref{tab:convergence_even_opt_upwind}, respectively. 
\revTwo{The convergence behavior of the other conserved quantities behave similarly and are, thus, not reported for the sake of brevity.}
For comparison, the same convergence analysis is performed using the traditional upwind SBP operators from \cite{Mattsson17}, with results reported in Tables~\ref{tab:convergence_odd_upwind} and \ref{tab:convergence_even_upwind}.

\begin{table}[H]
\centering
\caption{ $log(l^2 \text{ errors})$ and convergence rates for the isentropic vortex using odd order ($3^\text{rd}$--$9^\text{th}$) boundary optimized upwind SBP operators} 
\label{tab:convergence_odd_opt_upwind} 
\begin{tabular}{lllllllll}
\toprule
  $m$  & $\log{l^2}_{(3)}$ & $q$ & $\log{l^2}_{(5)}$& $q$& $\log{l^2}_{(7)}$& $q$ & $\log{l^2}_{(9)}$& $q$ \\
\midrule
$ 34$ & -2.89 & --- & -3.26 & --- & -3.36 & --- & -3.45 & --- \\
$68$ & --3.73 & 2.78 & -4.57 & 4.33 & -4.87 & 5.03 & -5.00 & 5.17 \\
$136$ & -4.63 & 3.01 & -5.83 & 4.19 & -6.31 & 4.77 & -6.53 & 5.06 \\
$272$ & -5.53 & 2.99 & -7.14 & 4.36 & -7.99 & 5.57 & -8.51 & 6.60 \\
$544$ & -6.41 & 2.98 & -8.49 & 4.47 & -9.81 & 6.08 & -10.72 & 7.32 \\
 \bottomrule
\end{tabular}
\end{table}
\begin{table}[H]
\centering
\caption{ $log(l^2 \text{ errors})$ and convergence rates for the isentropic vortex using even order ($2^\text{nd}$--$8^\text{th}$) boundary optimized upwind SBP operators} 
\label{tab:convergence_even_opt_upwind} 
\begin{tabular}{lllllllll}
\toprule
  $m$  & $\log{l^2}_{(2)}$ & $q$ & $\log{l^2}_{(4)}$& $q$& $\log{l^2}_{(6)}$& $q$ & $\log{l^2}_{(8)}$& $q$ \\
\midrule
$ 34$ & -2.36 & --- & -3.13 & --- & -3.40 & --- & -3.51 & --- \\
$68$ & --2.85 & 1.61 & -4.00 & 2.92 & -4.49 & 3.62 & -4.71 & 3.99 \\
$136$ & -3.41 & 1.86 & -5.07 & 3.55 & -5.95 & 4.83 & -6.30 & 5.26 \\
$272$ & -4.00 & 1.98 & -6.27 & 3.99 & -7.61 & 5.53 & -8.19 & 6.29 \\
$544$ & -4.61 & 2.00 & -7.49 & 4.04 & -9.41 & 5.99 & -10.38 & 7.28 \\
 \bottomrule
\end{tabular}
\end{table}
\begin{table}[H]
\centering
\caption{ $log(l^2 \text{ errors})$ and convergence rates for the isentropic vortex using odd order ($3^\text{rd}$--$9^\text{th}$) upwind SBP operators \cite{Mattsson17}} 
\label{tab:convergence_odd_upwind} 
\begin{tabular}{lllllllll}
\toprule
  $m$  & $\log{l^2}_{(3)}$ & $q$ & $\log{l^2}_{(5)}$& $q$& $\log{l^2}_{(7)}$& $q$ & $\log{l^2}_{(9)}$& $q$ \\
\midrule
$ 34$ & -2.88 & --- & -3.22 & --- & -3.26 & --- & -3.13 & --- \\
$68$ & --3.69 & 2.69 & -4.40 & 3.92 & -4.53 & 4.25 & -4.44 & 4.36 \\
$136$ & -4.54 & 2.82 & -5.52 & 3.71 & -5.76 & 4.09 & -6.15 & 5.67 \\
$272$ & -5.36 & 2.73 & -6.58 & 3.50 & -7.02 & 4.17 & -7.94 & 5.95 \\
$544$ & -6.16 & 2.65 & -7.59 & 3.39 & -8.28 & 4.18 & -9.68 & 5.8 \\
 \bottomrule
\end{tabular}
\end{table}
\begin{table}[H]
\centering
\caption{ $log(l^2 \text{ errors})$ and convergence rates for the isentropic vortex using even order ($2^\text{nd}$--$8^\text{th}$) upwind SBP operators \cite{Mattsson17}} 
\label{tab:convergence_even_upwind} 
\begin{tabular}{lllllllll}
\toprule
  $m$  & $\log{l^2}_{(2)}$ & $q$ & $\log{l^2}_{(4)}$& $q$& $\log{l^2}_{(6)}$& $q$ & $\log{l^2}_{(8)}$& $q$ \\
\midrule
$ 34$ & -2.36 & --- & -3.09 & --- & -3.22 & --- & -3.13 & --- \\
$68$ & --2.85 & 1.62 & -3.99 & 2.99 & -4.40 & 3.93 & -4.42 & 4.31 \\
$136$ & -3.41 & 1.86 & -5.06 & 3.55 & -5.73 & 4.42 & -6.13 & 5.65 \\
$272$ & -4.00 & 1.98 & -6.24 & 3.91 & -7.04 & 4.37 & -7.93 & 6.00 \\
$544$ & -4.61 & 2.00 & -7.40 & 3.86 & -8.32 & 4.23 & -9.73 & 5.95 \\
 \bottomrule
\end{tabular}
\end{table}

The convergence studies presented in Tables~\ref{tab:convergence_odd_opt_upwind}--\ref{tab:convergence_even_upwind} demonstrate that the boundary-optimized upwind SBP operators achieve the expected convergence rates. For the odd-order boundary-optimized upwind SBP operators derived herein the convergence rate may exhibit slightly better than expected convergence behavior.. The two-dimensional nonlinear convergence tests reflect, somewhat, the enhanced convergence behavior observed for the one-dimensional linear results in Section~\ref{sec:Convergence}. 
\revTwo{In contrast to the linear wave propagation results, the improvement in the observed convergence rates and error magnitudes in Tables~\ref{tab:convergence_odd_opt_upwind}--\ref{tab:convergence_even_upwind} for this nonlinear problem is less pronounced at moderate resolutions. Nevertheless, at the highest resolution considered, a clear improvement is observed for the boundary-optimized operators. For example, consider the ninth-order operators with $m=544$. The boundary-optimized operator achieves an observed convergence rate of $7.32$, compared to $5.80$ for the corresponding traditional SBP operator. The respective errors are approximately $10^{-10.72}$ and $10^{-9.68}$, which means that the error reduces by a factor of $\approx11$. Similarly, for the eighth-order operators, the errors at $m=544$ are approximately $10^{-10.38}$ and $10^{-9.73}$ for the boundary-optimized and traditional SBP operators, respectively. This corresponds to a reduction by a factor of $\approx 4.5$. 
}

Next, the robustness and accuracy properties of the boundary-optimized upwind SBP operators are examined. Time integration is performed using error-based step-size control with absolute and relative tolerances set to $10^{-7}$. Smaller tolerances were also tested, but no qualitative differences in robustness were observed. For this comparison, an $8^{\text{th}}$-order diagonal-norm SBP operator is employed on a $61 \times 61$ grid in each block, such that the isentropic vortex is only marginally resolved.

In addition to the boundary-optimized upwind SBP operators developed in this work, the following operators are considered for comparison: traditional SBP operators \cite{MattssonNordstrom04}, boundary-optimized traditional SBP operators \cite{MattssonAlmquistWeide2018}, and the upwind SBP operators from \cite{Mattsson17}. For the initial comparison, the same vortex parameters are used, and the solution is advanced to a final time of $t = 75$, corresponding to approximately seven vortex passages through the periodic domain.

\begin{figure}[H]
\centering

\begin{subfigure}[b]{0.47\textwidth}
  \centering
  \includegraphics[width=\linewidth, trim=0 0 0 0, clip]{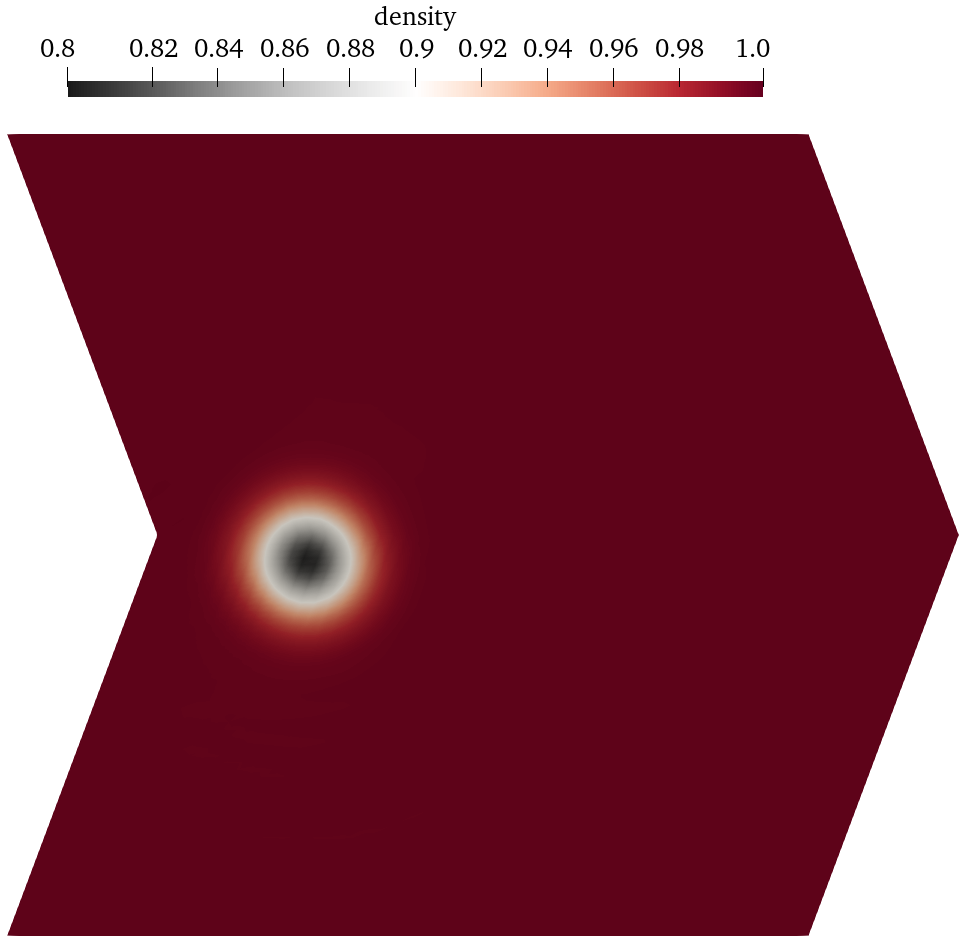}
  \caption{Boundary-optimized upwind SBP operator}
\end{subfigure}
\hspace*{0.5cm}
\begin{subfigure}[b]{0.47\textwidth}
  \centering
  \includegraphics[width=\linewidth, trim=0 0 0 0, clip]{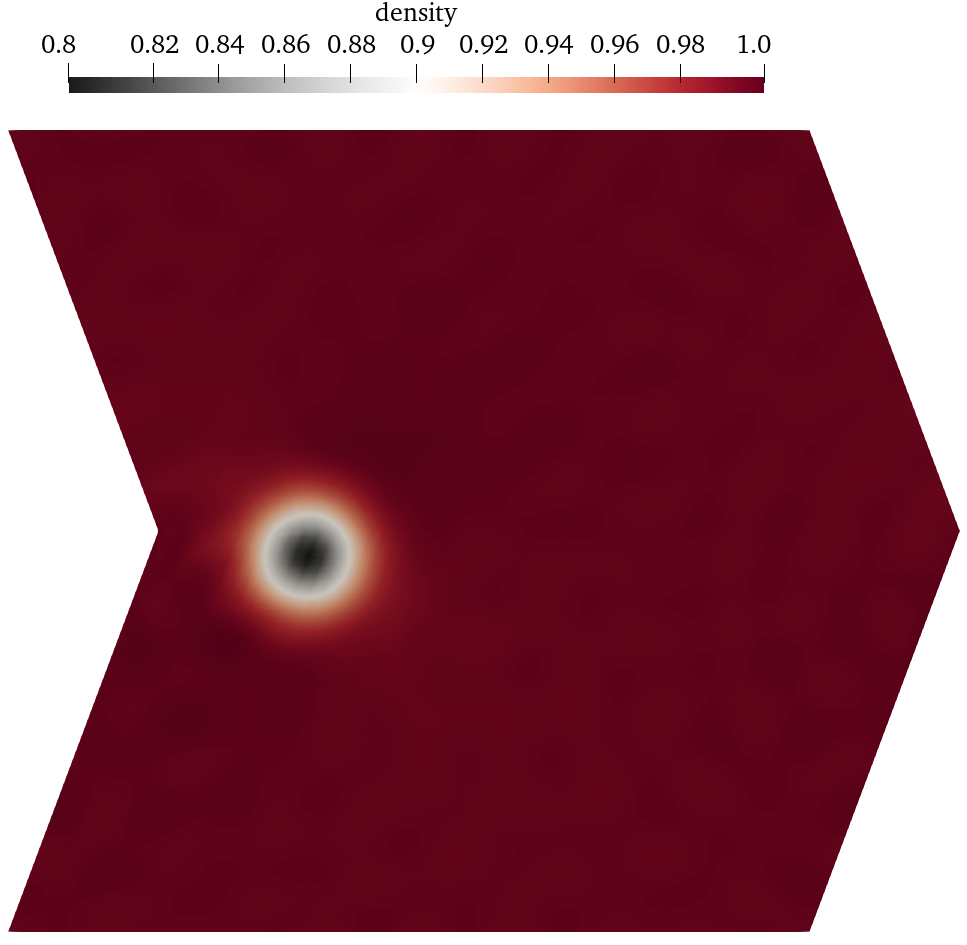}
  \caption{Traditional upwind SBP operator \cite{Mattsson17}}
\end{subfigure}
\begin{subfigure}[b]{0.47\textwidth}
  \includegraphics[width=\linewidth, trim=0 0 0 0, clip]{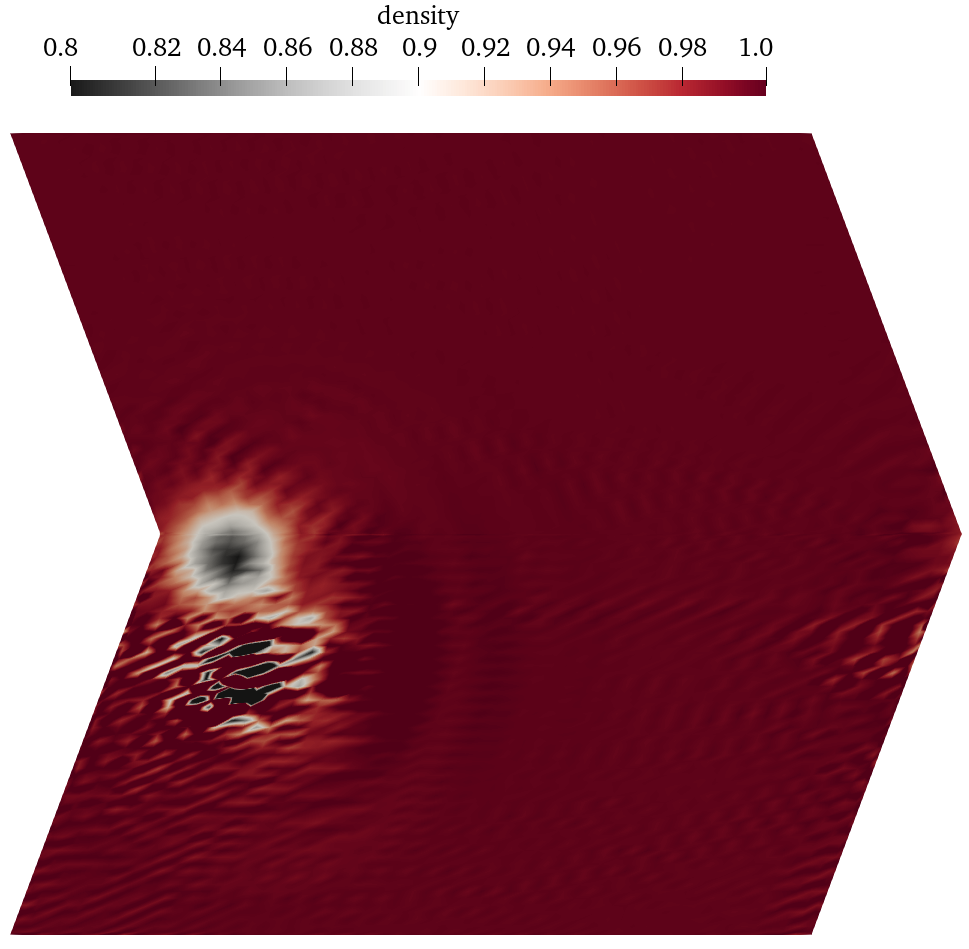}
  \caption{Boundary optimized traditional SBP operator \cite{MattssonAlmquistWeide2018}}
\end{subfigure}
\hspace*{0.5cm}
\begin{subfigure}[b]{0.47\textwidth}
  \centering
  \includegraphics[width=\linewidth, trim=0 0 0 0, clip]{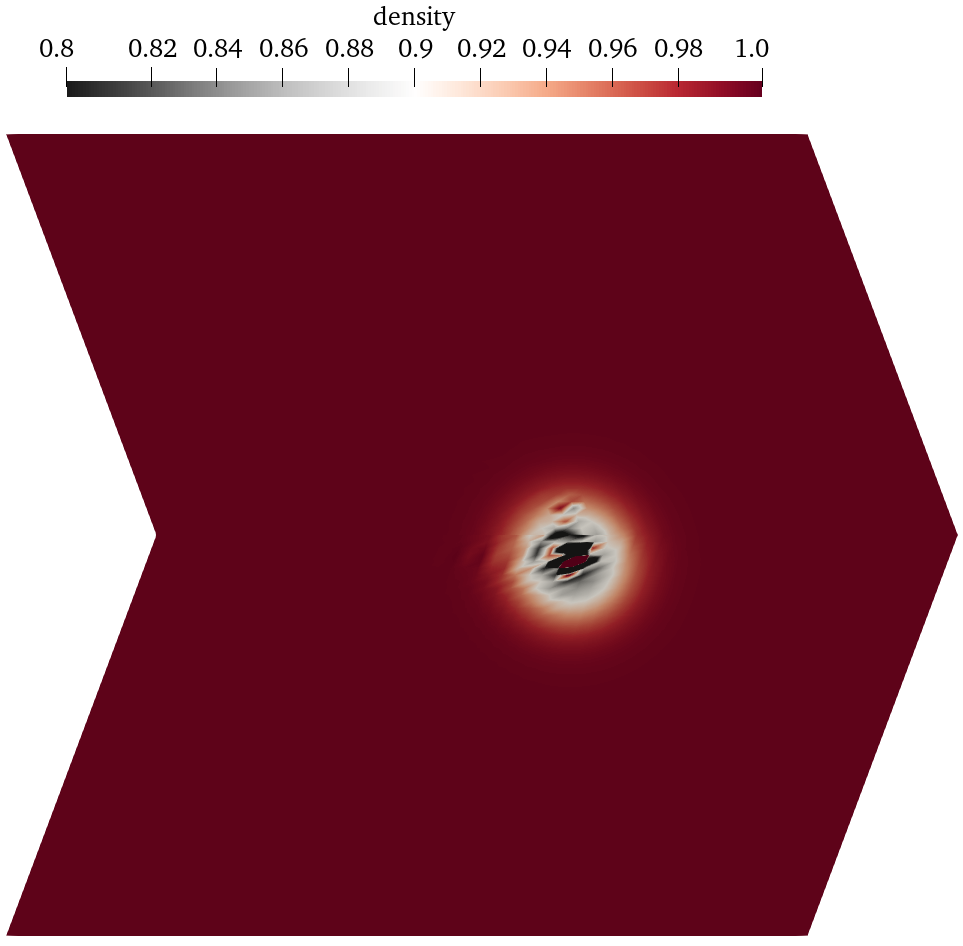}
  \caption{Traditional SBP operator \cite{MattssonNordstrom04}}
\end{subfigure}

\caption{
  Isentropic vortex on a chevron domain with two blocks each containing $61\times 61$ nodes. The vortex strength is $\varepsilon=5$. All results use $8^{\text{th}}$ order diagonal norm SBP operators with $4^{\text{th}}$ order boundary stencils. The boundary optimized upwind (a) and traditional upwind SBP operators run to the final time $t=75$, although the tranditional upwind operator exhibits small discrepancies. The boundary optimized traditional SBP operator crashes at $t\approx 41.438$ and the traditional SBP operator crashes early in the simulation at $t\approx 0.2919$.}
\label{fig:crash_test1}
\end{figure}

Figure~\ref{fig:crash_test1} illustrates the vortex solution obtained using each operator. \revTwo{When we report instability (and thus failure) of a particular run, it means that the density became negative during the approximation.} The traditional SBP operator becomes unstable early in the simulation at $t \approx 0.2919$, whereas the boundary-optimized traditional SBP operator remains stable until approximately $t \approx 41.438$. This indicates that improved boundary accuracy enhances the robustness of the discretization. Both the traditional and boundary-optimized upwind SBP operators remain stable until the final time. However, small discrepancies are observed in the solution obtained with the traditional upwind operator that are not present when using the boundary-optimized upwind SBP operator.

\begin{figure}[H]
\centering

\begin{subfigure}[b]{0.47\textwidth}
  \centering
  \includegraphics[width=\linewidth, trim=0 0 0 0, clip]{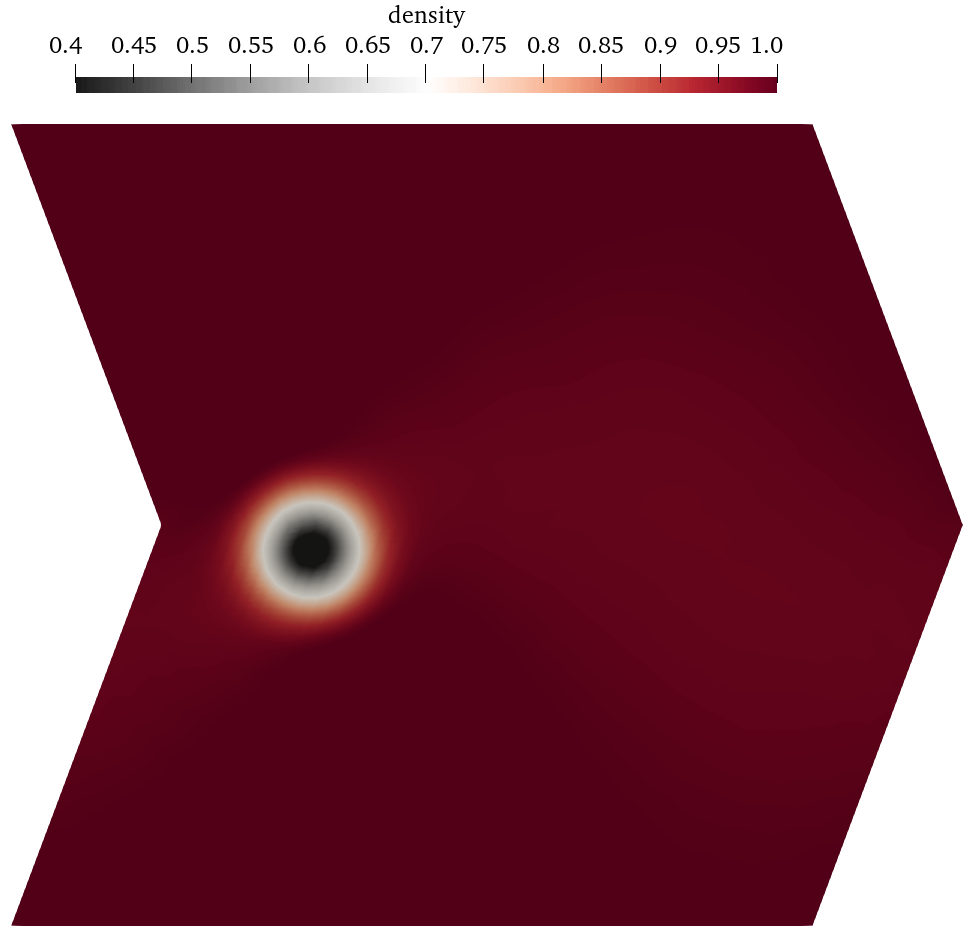}
  \caption{Boundary-optimized upwind SBP operator}
\end{subfigure}
\hspace*{0.5cm}
\begin{subfigure}[b]{0.47\textwidth}
  \centering
  \includegraphics[width=\linewidth, trim=0 0 0 0, clip]{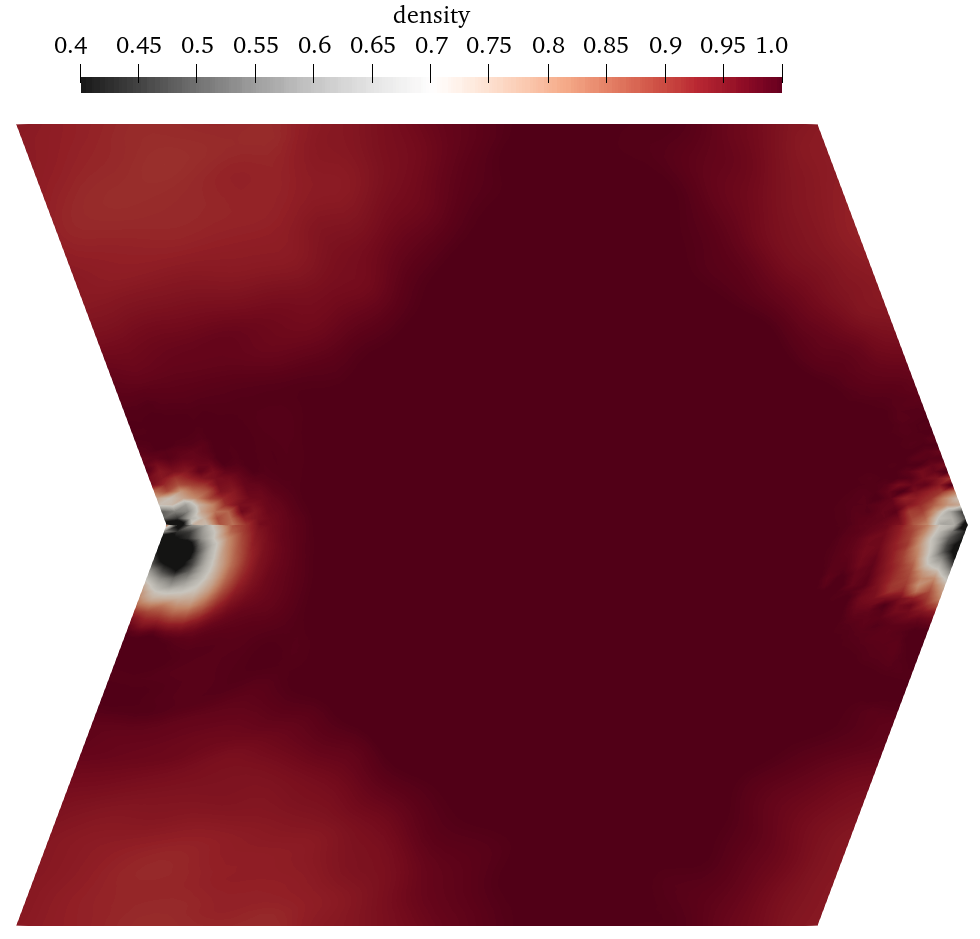}
  \caption{Traditional upwind SBP operator \cite{Mattsson17}}
\end{subfigure}

\caption{
  Isentropic vortex on a chevron domain with two blocks each containing $61\times 61$ nodes. The vortex strength is $\varepsilon=10$. The $8^{\text{th}}$ order boundary optimized upwind SBP operator (left) runs to the final time $t=75$ whereas the traditional $8^{\text{th}}$ order upwind SBP operator (right) crashes at $t\approx 72.255$.
}
\label{fig:crash_test2}
\end{figure}

To further investigate this behavior, the vortex strength is increased by setting $\varepsilon = 10$. The simulations are again performed using both the traditional and boundary-optimized upwind SBP operators and advanced to the final time $t = 75$. The boundary-optimized operators remain stable throughout the simulation, whereas the traditional upwind SBP operator becomes unstable at approximately $t \approx 72.255$, as shown in Figure~\ref{fig:crash_test2}.

As in the previous experiment, the instability originates near the block interface where the vortex intersects the boundary. This observation further indicates that improved boundary accuracy enhances the robustness of the discretization.

\subsection{Kelvin-Helmholtz instability}

As a final test, we consider a Kelvin--Helmholtz instability (KHI) setup for the two-dimensional compressible Euler equations of an ideal fluid to further examine the robustness of the boundary-optimized upwind SBP operators for under-resolved flows. Specifically, we employ the initial condition proposed by Chan et al.~\cite{chan2022entropy},
\begin{equation}
  \rho = \frac{1}{2} + \frac{3}{4} B(x,y),
  \quad
  v_1 = \frac{1}{2} \bigl( B(x,y) - 1 \bigr),
  \quad
  v_2 = \frac{1}{10}\sin(2 \pi x),
  \quad
  p = 1,
\end{equation}
where $B(x,y)$ is the smoothed approximation
\begin{equation}
  B(x, y) = \tanh(15 y + 7.5) - \tanh(15 y - 7.5)
\end{equation}
to a discontinuous step function. The computational domain is $[-1,1]^2$ with time interval $t \in [0,15]$.

Among the available flux vector splittings for the compressible Euler equations, the Steger--Warming splitting \cite{steger1979flux} is employed for the KHI tests. Alternative splittings, including the local Lax--Friedrichs and van Leer--H\"anel splittings \cite{vanleer1982flux,hanel1987accuracy,liou1991high}, were also tested and produced qualitatively similar results.

Time integration of the semidiscrete system is performed using the third-order, four-stage strong-stability-preserving (SSP) Runge--Kutta method \cite{kraaijevanger1991contractivity, conde2022embedded} combined with the error-based step-size controller proposed in \cite{ranocha2021optimized}. Error-based time stepping is applied with absolute and relative tolerances set to $10^{-6}$.

A resolution study similar to those presented in~\cite{ranocha2025robustness, glaubitz2025generalized} is conducted. The number of nodes per element is fixed while the number of elements is increased to refine the resolution. Specifically, the number of nodes is set to $m=17$ in each spatial direction, corresponding to the minimum number of points required to construct boundary-optimized upwind SBP operators with interior orders $p=8$ and $p=9$.

Table~\ref{tab:kelvin_helmholtz} summarizes the final simulation times obtained using the boundary-optimized upwind SBP operators and the traditional upwind SBP operators from \cite{Mattsson17}.

\begin{table}[htbp]
\centering
  \caption{Final times of numerical simulations of the Kelvin-Helmholtz
           instability with $K$ elements using $17$ nodes per coordinate
           direction for the upwind SBP methods. Final times less than 15
           indicate that the simulation crashed.}
  \label{tab:kelvin_helmholtz}
  \begin{subtable}{0.45\textwidth}
  \centering
    \caption{Boundary-optimized upwind SBP operator}
    \footnotesize
    \begin{tabular}{r rrrrrrrr}
      \toprule
      \multicolumn{1}{c}{$K$} & \multicolumn{8}{c}{interior order of accuracy} \\
          & \multicolumn{1}{c}{2}
          & \multicolumn{1}{c}{3}
          & \multicolumn{1}{c}{4}
          & \multicolumn{1}{c}{5}
          & \multicolumn{1}{c}{6}
          & \multicolumn{1}{c}{7} 
          & \multicolumn{1}{c}{8} 
          & \multicolumn{1}{c}{9} \\
      \midrule
        1  & 15 & 15 & 15 & 15 & 15 & 15 & 15 & 15 \\
        4  & 15 & 15 & 15 & 15 & 15 & 15 & 15 & 15 \\
       16  & 15 & 15 & 15 & 15 & 15 & 15 &  4.06 &  4.95 \\
       64  & 15 & 15 & 15 &  3.83 &  3.89 &  3.84 &  4.27 &  4.23 \\
      256  & 15 & 15 & 15 &  3.70 & 15 &  3.60 &  3.63 &  3.56 \\
     \bottomrule
    \end{tabular}
  \end{subtable}%
  \hspace*{0.5cm}
  \begin{subtable}{0.45\textwidth}
  \centering
    \caption{Traditional upwind SBP operator \cite{Mattsson17}}
    \footnotesize
    \begin{tabular}{r rrrrrrrr}
      \toprule
      \multicolumn{1}{c}{$K$} & \multicolumn{8}{c}{interior order of accuracy} \\
          & \multicolumn{1}{c}{2}
          & \multicolumn{1}{c}{3}
          & \multicolumn{1}{c}{4}
          & \multicolumn{1}{c}{5}
          & \multicolumn{1}{c}{6}
          & \multicolumn{1}{c}{7} 
          & \multicolumn{1}{c}{8} 
          & \multicolumn{1}{c}{9} \\
      \midrule
        1  & 15 & 15 & 15 & 15 & 15 & 15 & 15 & 15 \\
        4  & 15 & 15 & 15 & 15 & 15 & 15 & 15 & 15 \\
       16  & 15 & 15 & 15 & 15 & 15 &  1.99 &  1.79 &  1.71 \\
       64  & 15 & 15 &  4.52 &  3.83 &  3.73 &  3.71 &  2.38 &  1.92 \\
      256  & 15 & 15 &  5.76 &  3.67 &  3.63 &  3.62 &  3.50 &  3.50 \\
     \bottomrule
    \end{tabular}
  \end{subtable}%
\end{table}
First, it is observed that all upwind SBP methods with a small number of elements, $K \in \{1,4\}$, successfully complete the simulation. Furthermore, the boundary-optimized upwind SBP operators derived herein reach the final simulation time for lower interior orders of accuracy, $p \in \{2,3,4\}$. Overall, the boundary-optimized upwind SBP operators demonstrate improved robustness compared with their traditional counterparts. However, consistent with the observations reported in \cite{ranocha2025robustness}, simulations employing larger numbers of elements and higher orders of accuracy become unstable and terminate before $t=5$ for all upwind SBP operators.

The improved robustness of the boundary-optimized upwind SBP operators can be attributed to increased accuracy near element boundaries. Instabilities in traditional upwind SBP operators typically originate near element interfaces, where unphysical solution states (e.g., negative densities) may arise and ultimately lead to simulation failure. The improved boundary accuracy of the boundary-optimized operators enhances the propagation of wave-like solution structures across element interfaces.

\begin{figure}[H]
\centering

\begin{subfigure}[b]{0.47\textwidth}
  \centering
  \includegraphics[width=\linewidth, trim=80 20 250 30, clip]{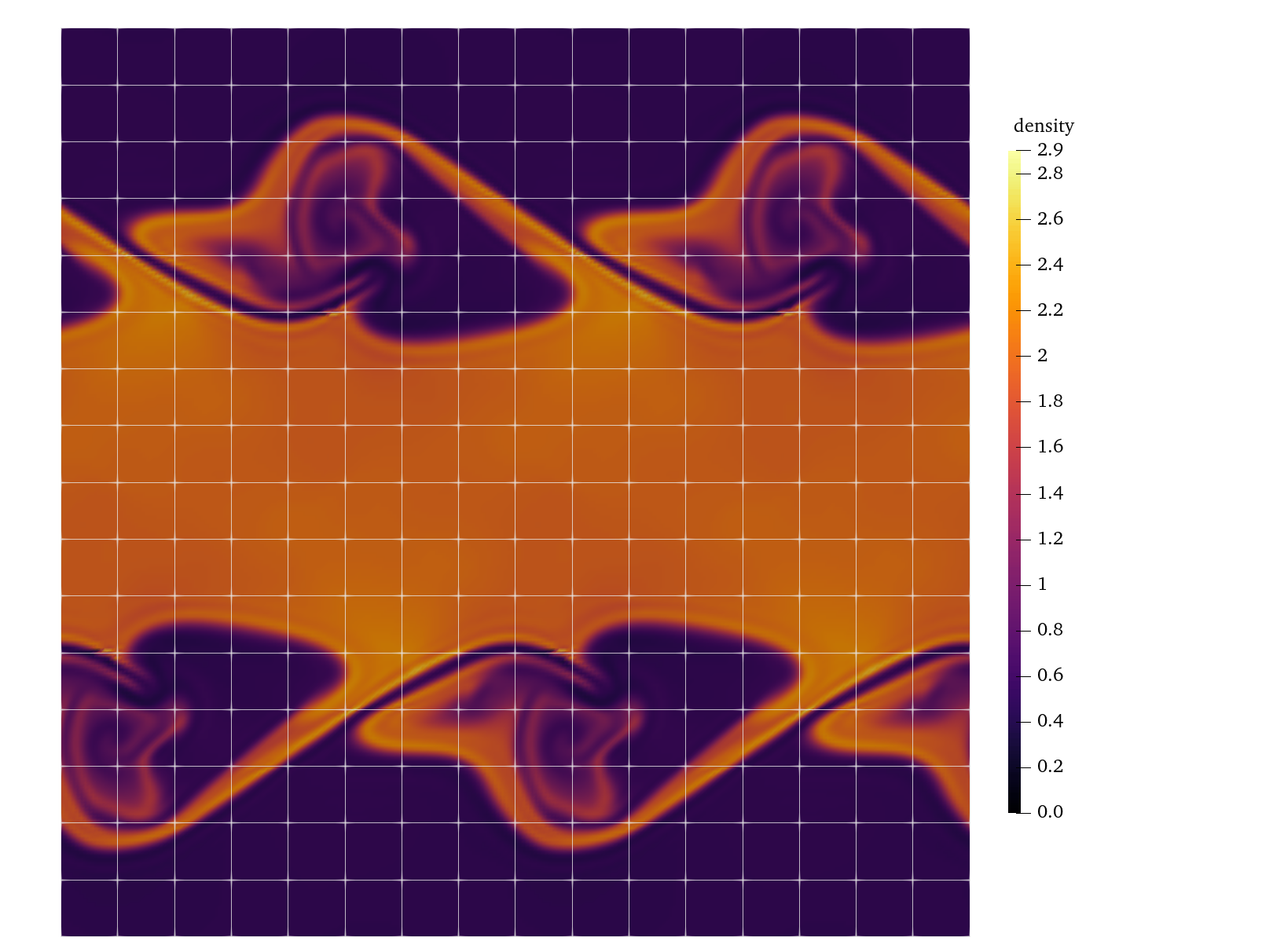}
  \caption{Boundary-optimized upwind SBP operator}
\end{subfigure}
\hspace*{0.5cm}
\begin{subfigure}[b]{0.47\textwidth}
  \centering
  \includegraphics[width=\linewidth, trim=80 20 250 30, clip]{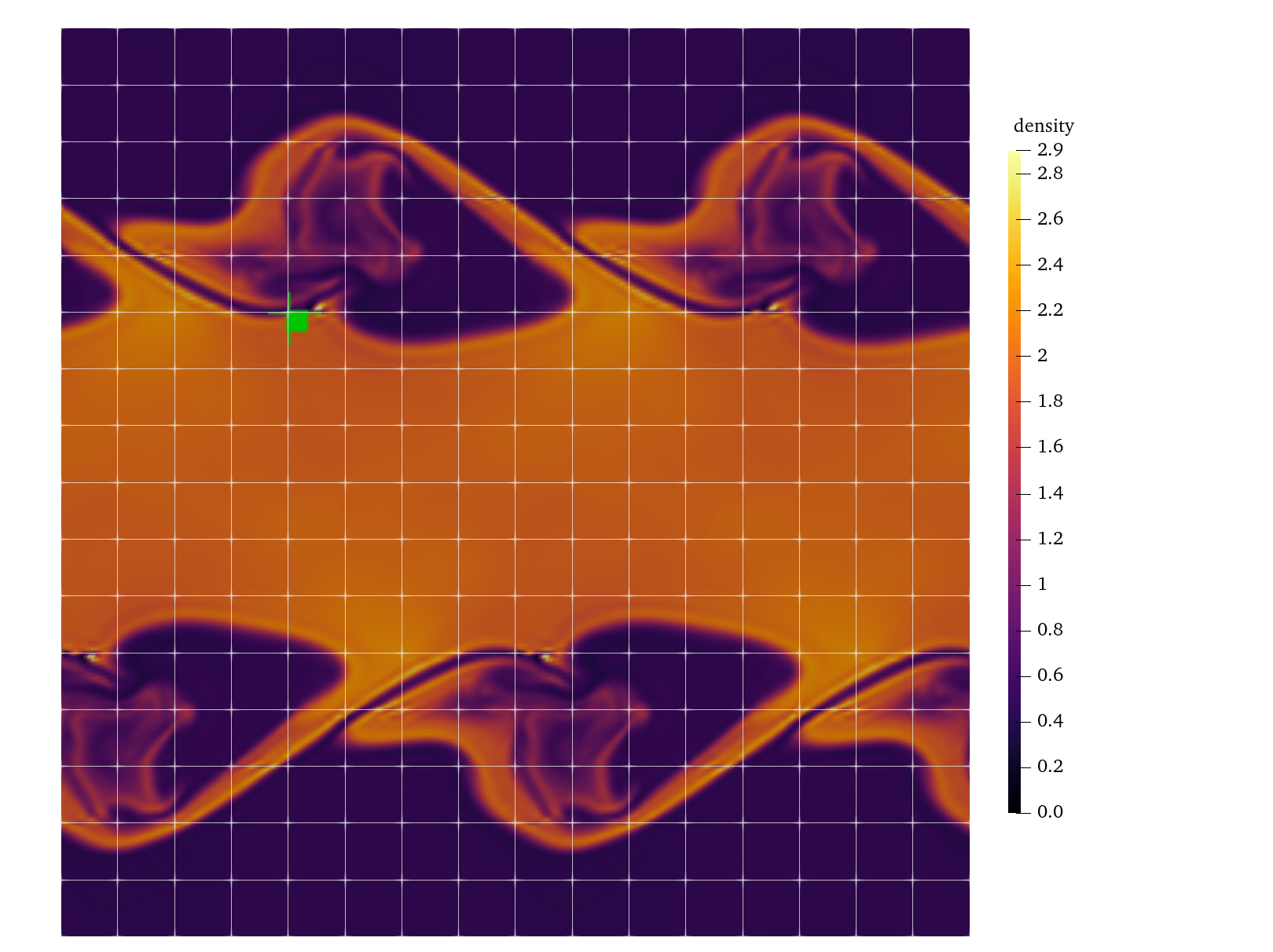}
  \caption{Traditional upwind SBP operator \cite{Mattsson17}}
\end{subfigure}

\caption{
  Visualization of numerical solutions from Kelvin--Helmholtz instability simulations at time $t = 3.63$.
  We use interior order $p=6$ upwind SBP operators on a mesh of $16^2$ Cartesian elements with $17^2$ nodes per element.
  The boundary-optimized simulation reaches the final time, whereas the traditional one crashes; green regions in (b) mark points with negative density.
}
\label{fig:kelvin_helmholtz}
\end{figure}

This behavior is illustrated in Figure~\ref{fig:kelvin_helmholtz}, which shows the Kelvin--Helmholtz solution at final time $t = 3.63$. The simulation uses a $16 \times 16$ Cartesian mesh (so $K = 256$) with $17$ nodes in each spatial direction per element and $6^{\textrm{th}}$ order accurate upwind SBP operators. From the last row in Table~\ref{tab:kelvin_helmholtz}, this configuration remains stable for the boundary-optimized upwind SBP operators but becomes unstable for the traditional upwind SBP operators. The regions highlighted in bright green in Figure~\ref{fig:kelvin_helmholtz}(b) correspond to locations where the density becomes negative in the simulation using the traditional upwind SBP operators. These regions are concentrated near element boundaries.

\section{Conclusions}\label{sec:Conclusions}

The primary objective of this work has been the construction of boundary-optimized diagonal-norm upwind SBP operators that provide built-in artificial dissipation when combined with flux vector splitting techniques for hyperbolic systems. Boundary and interface conditions are imposed using either the projection method or the simultaneous approximation term (SAT) technique.

The proposed boundary-optimized upwind SBP operators yield robust and accurate discretizations for hyperbolic problems, as demonstrated through numerical experiments involving both linear and nonlinear systems. For wave propagation problems, the numerical results show that the boundary-optimized diagonal-norm upwind SBP operators outperform both boundary-optimized central SBP operators and traditional upwind SBP operators. The observed convergence rates for linear first-order wave equations exceed the theoretical expectations. \revTwo{These observed convergence rates that are higher than the rates predicted from the formal boundary accuracy of the operators over the resolutions considered, are attributed to suppression of the leading boundary truncation-error terms.}

Improved convergence behavior and accuracy, particularly for marginally resolved discretizations, are also observed for nonlinear two-dimensional simulations of the compressible Euler equations. Furthermore, the boundary-optimized upwind SBP operators demonstrate enhanced robustness for complex but smooth Kelvin--Helmholtz instability simulations compared with traditional upwind SBP operators based on uniform-grid stencils.

Future work will investigate the efficiency and robustness of the proposed upwind SBP operators in more complex settings, including curvilinear geometries and three-dimensional domains.

\section*{CRediT authorship contribution statement}
	{\bf Ken Mattsson:} Conceptualization; Formal analysis; Investigation; Methodology; Data curation; Visualization; Software; Writing - original draft
	{\bf David Niemel\"{a}:} Conceptualization; Methodology; Software; Writing - review \& editing
	{\bf Andrew R. Winters:} Conceptualization; Visualization; Software; Writing - review \& editing
	
\section*{Data Availability}

	Data is available upon request to the corresponding author.
	
\section*{Declaration of competing interest}
	The authors declare that they have no known competing financial interests or personal relationships that could have appeared to influence the work reported in this paper.
	
\section*{Acknowledgement}	
\noindent	
	Ken Mattsson was partially supported by the Swedish Research Council (grant 2021-05830 VR), and FORMAS (grant 2022-00843). \\
	David Niemel\"a was supported by the Swedish Research Council (grant 2021-05830 VR) and FORMAS (grant 2022-00843). \\
	Andrew Winters was supported by the Swedish Research Council (grant 2020-03642 VR).

\bibliographystyle{elsarticle-num}
\bibliography{referenser}

@string{Bertil 	= {B. Gustafsson}}

@string{Kreiss 	= {H--O Kreiss}}

@string{Mark  	= {M. H. Carpenter}}

@string{BIT 	= {BIT}}

@article{MattssonSvard04,
  author      = {K. Mattsson and M. Sv{\"a}rd and J. Nordstr{\"o}m},
  title       = {{S}table and {A}ccurate {A}rtificial {D}issipation},
  journal     = {Journal of Scientific Computing},
  month       = {August},
  volume      = {21(1)},
  year        = 2004,
    pages		=	{57-79},

}

@article{MattssonNordstrom04,
  author      = {K. Mattsson and J. Nordstr{\"o}m},
  title       = { Summation by parts operators for finite difference
                  approximations of second derivatives},
  journal     = {J. Comput. Phys.},
  year        = {2004},
  volume      = {199(2)},
  pages		=	{503-540},
}

@article{Mattsson11,
year={2012},
issn={0885-7474},
journal={Journal of Scientific Computing},
volume={51},
issue={3},
doi={10.1007/s10915-011-9525-z},
title={Summation by Parts Operators for Finite Difference Approximations of Second-Derivatives with Variable Coefficients},
url={http://dx.doi.org/10.1007/s10915-011-9525-z},
publisher={Springer US},
author={Mattsson, Ken},
pages={650-682},
language={English}
}

@article{MattssonAlmquist13,
  author      = {K. Mattsson and M. Almquist} ,
  title       = {A solution to the stability issues with block norm summation by parts operators },
  journal     = {J. Comput. Phys.},
  year        = {2013},
  volume      = {253},
  pages		=	{418-442},
}

@article{MattssonAlmquistCarpenter13,
  author      = {K. Mattsson and M. Almquist and M. H. Carpenter} ,
  title       = {Optimal diagonal-norm {SBP} operators },
  journal     = {J. Comput. Phys.},
  year        = {2014},
  volume      = {264},
  pages       = {91-111}
}

@article{Mattsson14,
	Author = {Ken Mattsson},
	Doi = {http://dx.doi.org/10.1016/j.jcp.2014.06.027},
	Issn = {0021-9991},
	Journal = {J. Comput. Phys.},
	Number = {0},
	Pages = {432 - 454},
	Title = {Diagonal-norm summation by parts operators for finite difference approximations of third and fourth derivatives},
	Url = {http://www.sciencedirect.com/science/article/pii/S0021999114004343},
	Volume = {274},
	Year = {2014}}

@article{Mattsson17,
title = "Diagonal-norm upwind {SBP} operators ",
journal = "Journal of Computational Physics ",
volume = "335",
number = "",
pages = "283 - 310",
year = "2017",
note = "",
issn = "0021-9991",
doi = "http://dx.doi.org/10.1016/j.jcp.2017.01.042",
url = "http://www.sciencedirect.com/science/article/pii/S002199911730058X",
author = "Ken Mattsson"
}

@article{MattssonOssian17,
title = "Compatible diagonal-norm staggered and upwind SBP operators",
journal = "Journal of Computational Physics",
volume = "352",
number = "Supplement C",
pages = "52 - 75",
year = "2018",
issn = "0021-9991",
doi = "https://doi.org/10.1016/j.jcp.2017.09.044",
url = "http://www.sciencedirect.com/science/article/pii/S0021999117307040",
author = "Ken Mattsson and Ossian O'Reilly"
}

@article{MattssonAlmquistWeide2018,
title = "Boundary optimized diagonal-norm {SBP} operators",
journal = "Journal of Computational Physics",
year = "2018",
issn = "0021-9991",
doi = "https://doi.org/10.1016/j.jcp.2018.06.010",
url = "http://www.sciencedirect.com/science/article/pii/S0021999118303887",
author = "Ken Mattsson and Martin Almquist and Edwin van der Weide"
}

@article{MattssonOlsson2018,
title = "An improved projection method",
journal = "Journal of Computational Physics",
volume = "372",
pages = "349 - 372",
year = "2018",
issn = "0021-9991",
doi = "https://doi.org/10.1016/j.jcp.2018.06.030",
url = "http://www.sciencedirect.com/science/article/pii/S002199911830408X",
author = "Ken Mattsson and Pelle Olsson"
}

@article{LundgrenMattsson_2020,
title = {An efficient finite difference method for the shallow water equations},
journal = {Journal of Computational Physics},
volume = {422},
pages = {109784},
year = {2020},
issn = {0021-9991},
doi = {https://doi.org/10.1016/j.jcp.2020.109784},
url = {https://www.sciencedirect.com/science/article/pii/S0021999120305581},
author = {Lukas Lundgren and Ken Mattsson}
}

@article{VidarMattsson_D2_2023,
title = {Boundary-optimized summation-by-parts operators for finite difference approximations of second derivatives with variable coefficients},
journal = {Journal of Computational Physics},
volume = {491},
pages = {112376},
year = {2023},
issn = {0021-9991},
doi = {https://doi.org/10.1016/j.jcp.2023.112376},
url = {https://www.sciencedirect.com/science/article/pii/S0021999123004710},
author = {Vidar Stiernstr{\"o}m and Martin Almquist and Ken Mattsson}
}

@article{VidarMattsson_RV_2021,
title = {A residual-based artificial viscosity finite difference method for scalar conservation laws},
journal = {Journal of Computational Physics},
volume = {430},
pages = {110100},
year = {2021},
issn = {0021-9991},
doi = {https://doi.org/10.1016/j.jcp.2020.110100},
url = {https://www.sciencedirect.com/science/article/pii/S0021999120308743},
author = {Vidar Stiernstr{\"o}m and Lukas Lundgren and Murtazo Nazarov and Ken Mattsson}
}

@article{ErikssonMattsson22,
	author = {Eriksson, Gustav and Mattsson, Ken},
	date = {2022/07/23},
	doi = {10.1007/s10915-022-01941-5},
	id = {Eriksson2022},
	isbn = {1573-7691},
	journal = {Journal of Scientific Computing},
	number = {3},
	pages = {81},
	title = {Weak Versus Strong Wall Boundary Conditions for the Incompressible Navier-Stokes Equations},
	url = {https://doi.org/10.1007/s10915-022-01941-5},
	volume = {92},
	year = {2022}}

@article{Svard2014,
title = "Review of summation-by-parts-operators schemes for initial-boundary-value problems ",
journal = "Journal of Computational Physics ",
volume = "268",
number = "0",
pages = "17 - 38",
year = "2014",
note = "",
issn = "0021-9991",
doi = "http://dx.doi.org/10.1016/j.jcp.2014.02.031",
url = "http://www.sciencedirect.com/science/article/pii/S002199911400151X",
author = {Magnus Sv\"ard and Jan Nordstr\"om}
}

@article{Svard04,
  author  =  {M. Sv\"ard},
  title   =  {On coordinate transformation for summation-by-parts operators },
  journal =  {Journal of Scientific Computing},
  volume  =   {20(1)},
  year    =   2004,
}

@article{SVARD2019108819,
title = {On the convergence rates of energy-stable finite-difference schemes},
journal = {Journal of Computational Physics},
volume = {397},
pages = {108819},
year = {2019},
issn = {0021-9991},
doi = {https://doi.org/10.1016/j.jcp.2019.07.018},
url = {https://www.sciencedirect.com/science/article/pii/S0021999119305030},
author = {Magnus Sv{\"a}rd and Jan Nordstr{\"o}m},
}

@article{DelReyFernandez2014a,
title = "Review of summation-by-parts operators with simultaneous approximation terms for the numerical solution of partial differential equations",
journal = "Computers \& Fluids",
volume = "95",
number = "",
pages = "171 - 196",
year = "2014",
note = "",
issn = "0045-7930",
doi = "http://dx.doi.org/10.1016/j.compfluid.2014.02.016",
url = "http://www.sciencedirect.com/science/article/pii/S0045793014000796",
author = "David C. Del Rey Fern\'{a}ndez and Jason E. Hicken and David W. Zingg"
}

@article{DelReyFernandez2014,
 author = {Del Rey Fern\'{a}ndez, David C. and Boom, Pieter D. and Zingg, David W.},
 title = {A Generalized Framework for Nodal First Derivative Summation-by-parts Operators},
 journal = {J. Comput. Phys.},
 issue_date = {June, 2014},
 volume = {266},
 month = jun,
 year = {2014},
 issn = {0021-9991},
 pages = {214--239},
 numpages = {26},
 url = {http://dx.doi.org/10.1016/j.jcp.2014.01.038},
 doi = {10.1016/j.jcp.2014.01.038},
 acmid = {2599287},
 publisher = {Academic Press Professional, Inc.},
 address = {San Diego, CA, USA},
}

@article{KreissOliger72,
 author  = {Heinz--Otto Kreiss and Joseph Oliger},
 title	 = {Comparison of accurate methods for the integration of
             hyperbolic equations},
 journal = {Tellus},
 volume  = {XXIV},
 pages   = {199-215},
 year	 = 1972,
 issue   = 3,
}

@incollection{KreissScherer74,
title = "Finite Element and Finite Difference Methods for Hyperbolic Partial Differential Equations ",
editor = "Boor, Carl de ",
booktitle = "Mathematical Aspects of Finite Elements in Partial Differential Equations ",
publisher = "Academic Press",
edition = "",
address = "",
year = "1974",
pages = "195 - 212",
isbn = "978-0-12-208350-1",
doi = "https://doi.org/10.1016/B978-0-12-208350-1.50012-1",
url = "https://www.sciencedirect.com/science/article/pii/B9780122083501500121",
author = "H.-O. Kreiss and G. Scherer"
}

@article{Strand94,
 author 	=	{B. Strand},
 title		=	{Summation by Parts for Finite Difference Approximations for $d/dx$},
 journal	=	{J. Comput. Physics},
 year		=	1994,
 volume  	=	110,
 pages		=	{47-67},
}

@book{Gustafsson2007,
author = {Gustafsson, Bertil},
doi = {10.1007/978-3-540-74993-6},
isbn = {978-3-540-74993-6},
title = {{High Order Difference Methods for Time Dependent PDE}},
url = {https://books.google.com/books?id=guaT5yPiLAgC},
year = {2007},
publisher={Springer}
}

@article{Linders1109585,
   author = {Linders, Viktor and Kupiainen, Marco and Nordstrom, Jan},
   institution = {Climate research - Rossby Centre},
   journal = {Journal of Computational Physics},
   pages = {160--176},
   title = {Summation-by-Parts operators with minimal dispersion error for coarse grid flow calculations},
   volume = {340},
   DOI = {10.1016/j.jcp.2017.03.039},
   year = {2017}
}

@article{duru2024,
  title={A dual-pairing summation-by-parts finite difference framework for nonlinear conservation laws},
  author={Duru, Kenneth and Stewart, Dougal and Lee, Nathan},
  journal={arXiv preprint arXiv:2411.06629},
  year={2024}
}

@article{hew2025strongly,
  title={Strongly stable dual-pairing summation by parts finite difference schemes for the vector invariant nonlinear shallow water equations--{I}: {N}umerical scheme and validation on the plane},
  author={Hew, Justin Kin Jun and Duru, Kenneth and Roberts, Stephen and Zoppou, Christopher and Ricardo, Kieran},
  journal={Journal of Computational Physics},
  volume={523},
  pages={113624},
  year={2025},
  publisher={Elsevier}
}

@article{williams2024full,
  title={Full-spectrum dispersion relation preserving summation-by-parts operators},
  author={Williams, Christopher and Duru, Kenneth},
  journal={SIAM Journal on Numerical Analysis},
  volume={62},
  number={4},
  pages={1565--1588},
  year={2024},
  publisher={SIAM}
}

@article{ranocha2025robustness,
  title={On the robustness of high-order upwind summation-by-parts methods for nonlinear conservation laws},
  author={Ranocha, Hendrik and Winters, Andrew R and Schlottke-Lakemper, Michael and {\"O}ffner, Philipp and Glaubitz, Jan and Gassner, Gregor J},
  journal={Journal of Computational Physics},
  volume={520},
  pages={113471},
  year={2025},
  publisher={Elsevier}
}

@article{glaubitz2025generalized,
  title={Generalized upwind summation-by-parts operators and their application to nodal discontinuous {G}alerkin methods},
  author={Glaubitz, Jan and Ranocha, Hendrik and Winters, Andrew R and Schlottke-Lakemper, Michael and {\"O}ffner, Philipp and Gassner, Gregor},
  journal={Journal of Computational Physics},
  volume={529},
  pages={113841},
  year={2025},
  publisher={Elsevier}
}

@article{ranocha2022adaptive,
  title={Adaptive numerical simulations with {T}rixi.jl:
         {A} case study of {J}ulia for scientific computing},
  author={Ranocha, Hendrik and Schlottke-Lakemper, Michael and Winters, Andrew Ross
          and Faulhaber, Erik and Chan, Jesse and Gassner, Gregor},
  journal={Proceedings of the JuliaCon Conferences},
  volume={1},
  number={1},
  pages={77},
  year={2022},
  doi={10.21105/jcon.00077},
  eprint={2108.06476},
  eprinttype={arXiv},
  eprintclass={cs.MS}
}

@article{schlottkelakemper2021purely,
  title={A purely hyperbolic discontinuous {G}alerkin approach for
         self-gravitating gas dynamics},
  author={Schlottke-Lakemper, Michael and Winters, Andrew R and
          Ranocha, Hendrik and Gassner, Gregor J},
  journal={Journal of Computational Physics},
  pages={110467},
  year={2021},
  month={06},
  volume={442},
  publisher={Elsevier},
  doi={10.1016/j.jcp.2021.110467},
  eprint={2008.10593},
  eprinttype={arXiv},
  eprintclass={math.NA}
}

@incollection{ahrens2005paraview,
  title={{ParaView}: {A}n End-User Tool for Large-Data Visualization},
  author={Ahrens, James and Geveci, Berk and Law, Charles},
  booktitle={The Visualization Handbook},
  publisher={Elsevier},
  year={2005},
  pages={717--731}
}

@article{ranocha2021sbp,
  title={{SummationByPartsOperators.jl}: {A} {J}ulia library of provably stable
         semidiscretization techniques with mimetic properties},
  author={Ranocha, Hendrik},
  journal={Journal of Open Source Software},
  year={2021},
  month={08},
  doi={10.21105/joss.03454},
  volume={6},
  number={64},
  pages={3454},
  publisher={The Open Journal},
  url={https://github.com/ranocha/SummationByPartsOperators.jl}
}

@techreport{steger1979flux,
  title={Flux Vector Splitting of the Inviscid Gasdynamic Equations With
         Application to Finite Difference Methods},
  author={Steger, Joseph L and Warming, RF},
  institution={NASA},
  year={1979},
  type={Technical Memorandum},
  number={NASA/TM-78605},
  address={NASA Ames Research Center; Moffett Field, CA, United States}
}

@techreport{shu1997essentially,
  title={Essentially Non-Oscillatory and Weighted Essentially Non-Oscillatory
         Schemes for Hyperbolic Conservation Laws},
  author={Shu, Chi-Wang},
  institution= {NASA},
  month={11},
  year={1997},
  type={Final Report},
  number={NASA/CR-97-206253},
  address={Institute for Computer Applications in Science and Engineering,
           NASA Langley Research Center, Hampton VA United States}
}

@inproceedings{vanleer1982flux,
  title={Flux-vector splitting for the {E}uler equations},
  author={van Leer, Bram},
  pages={507--512},
  year={1982},
  booktitle={Eighth International Conference on Numerical Methods in Fluid
             Dynamics},
  editor={Krause, Egon},
  series={Lecture Notes in Physics},
  volume={170},
  publisher={Springer},
  address={Berlin, Heidelberg},
  doi={10.1007/3-540-11948-5_66}
}

@inproceedings{hanel1987accuracy,
  title={On the accuracy of upwind schemes for the solution of the
         {N}avier-{S}tokes equations},
  author={H{\"a}nel, D and Schwane, R and Seider, G},
  booktitle={8th Computational Fluid Dynamics Conference},
  pages={1105},
  year={1987},
  organization={American Institute of Aeronautics and Astronautics},
  doi={10.2514/6.1987-1105}
}

@techreport{liou1991high,
  title={High-Order Polynomial Expansions ({HOPE}) for flux-vector splitting},
  author={Liou, Meng-Sing and Steffen Jr, Chris J},
  institution={NASA},
  year={1991},
  type={Technical Memorandum},
  number={NASA/TM-104452},
  address={NASA Lewis Research Center; Cleveland, OH, United States}
}

@article{rackauckas2017differentialequations,
  title={{DifferentialEquations.jl} {--} {A} Performant and Feature-Rich
         Ecosystem for Solving Differential Equations in {J}ulia},
  author={Rackauckas, Christopher and Nie, Qing},
  journal={Journal of Open Research Software},
  volume={5},
  number={1},
  pages={15},
  year={2017},
  publisher={Ubiquity Press},
  doi={10.5334/jors.151}
}

@article{chan2022entropy,
  title={On the entropy projection and the robustness of high order entropy
         stable discontinuous {G}alerkin schemes for under-resolved flows},
  author={Chan, Jesse and Ranocha, Hendrik and Rueda-Ramirez, Andr{\'e}s M
          and Gassner, Gregor J and Warburton, Tim},
  journal={Frontiers in Physics},
  year={2022},
  month={07},
  doi={10.3389/fphy.2022.898028},
  volume={10},
  eprint={2203.10238},
  eprinttype={arxiv},
  eprintclass={math.NA}
}

@article{kraaijevanger1991contractivity,
  title={Contractivity of {R}unge-{K}utta methods},
  author={Kraaijevanger, Johannes Franciscus Bernardus Maria},
  journal={BIT Numerical Mathematics},
  volume={31},
  number={3},
  pages={482--528},
  year={1991},
  publisher={Springer},
  doi={10.1007/BF01933264}
}

@article{conde2022embedded,
  title={Embedded pairs for optimal explicit strong stability preserving
         {R}unge-{K}utta methods},
  author={Conde, Sidafa and Fekete, Imre and Shadid, John N},
  journal={Journal of Computational and Applied Mathematics},
  volume={412},
  pages={114325},
  year={2022},
  doi={10.1016/j.cam.2022.114325},
  eprint={1806.08693},
  eprinttype={arxiv},
  eprintclass={math.NA}
}

@article{ranocha2021optimized,
  title={Optimized {R}unge-{K}utta Methods with Automatic Step Size Control
         for Compressible Computational Fluid Dynamics},
  author={Ranocha, Hendrik and Dalcin, Lisandro and Parsani, Matteo
          and Ketcheson, David I.},
  journal={Communications on Applied Mathematics and Computation},
  volume={4},
  pages={1191--1228},
  year={2021},
  month={11},
  doi={10.1007/s42967-021-00159-w},
  eprint={2104.06836},
  eprinttype={arxiv},
  eprintclass={math.NA}
}

@article{ham2006accurate,
  title={Accurate and stable finite volume operators for unstructured flow solvers},
  author={Ham, Frank and Mattsson, Ken and Iaccarino, Gianluca},
  journal={Annual Research Briefs},
  volume={243},
  year={2006},
  publisher={Center for Turbulence Research, NASA Ames/Stanford University Stanford}
}

@article{nordstrom2003finite,
  title={Finite volume methods, unstructured meshes and strict stability for hyperbolic problems},
  author={Nordstr{\"o}m, Jan and Forsberg, Karl and Adamsson, Carl and Eliasson, Peter},
  journal={Applied Numerical Mathematics},
  volume={45},
  number={4},
  pages={453--473},
  year={2003},
  publisher={Elsevier}
}

@article{gassner2013skew,
  title={A skew-symmetric discontinuous {G}alerkin spectral element discretization and its relation to {SBP}-{SAT} finite difference methods},
  author={Gassner, Gregor J},
  journal={SIAM Journal on Scientific Computing},
  volume={35},
  number={3},
  pages={A1233--A1253},
  year={2013},
  publisher={SIAM}
}

@incollection{winters2021,
author={Winters, Andrew R. and Kopriva, David A. and Gassner, Gregor J. and Hindenlang, Florian},
editor={Kronbichler, Martin and Persson, Per-Olof},
title={Construction of Modern Robust Nodal Discontinuous {G}alerkin Spectral Element Methods for the Compressible {N}avier--{S}tokes Equations},
bookTitle={Efficient High-Order Discretizations for Computational Fluid Dynamics},
year={2021},
publisher={Springer International Publishing},
pages={117--196}
}

@book{butcher2008numerical,
  author    = {Butcher, J.C.},
  title     = {Numerical Methods for Ordinary Differential Equations},
  publisher = {John Wiley \& Sons},
  year      = {2008},
  edition   = {2nd}
}

@article{abgrall2020analysis,
	author = {Abgrall, R{\'e}mi and Nordstr{\"o}m, Jan and {\"O}ffner, Philipp and Tokareva, Svetlana},
	journal = {Journal of Scientific Computing},
	number = {2},
	pages = {1--29},
	publisher = {Springer},
	title = {Analysis of the {SBP-SAT} stabilization for finite element methods part {I}: linear problems},
	volume = {85},
	year = {2020}}

@article{Hanifa2025,
	author = {Hanif, Hanifa and Nordstr{\"o}m, Jan and Malan, Arnaud G.},
	journal = {AIMS Mathematics},
	number = {9},
	pages = {22579-22597},
	title = {Efficiency analysis of continuous and discontinuous {G}alerkin finite element methods},
	volume = {10},
	year = {2025}}

@article{del2019extension,
  title={Extension of tensor-product generalized and dense-norm summation-by-parts operators to curvilinear coordinates},
  author={Del Rey Fern{\'a}ndez, David C and Boom, Pieter D and Carpenter, Mark H and Zingg, David W},
  journal={Journal of Scientific Computing},
  volume={80},
  number={3},
  pages={1957--1996},
  year={2019},
  publisher={Springer}
}

@article{chan2018discretely,
  title={On discretely entropy conservative and entropy stable discontinuous {G}alerkin methods},
  author={Chan, Jesse},
  journal={Journal of Computational Physics},
  volume={362},
  pages={346--374},
  year={2018},
  publisher={Elsevier}
}

@article{ranocha2017extended,
  title={Extended skew-symmetric form for summation-by-parts operators and varying {J}acobians},
  author={Ranocha, Hendrik and {\"O}ffner, Philipp and Sonar, Thomas},
  journal={Journal of Computational Physics},
  volume={342},
  pages={13--28},
  year={2017},
  publisher={Elsevier}
}

@article{hicken2013summation,
  title={Summation-by-parts operators and high-order quadrature},
  author={Hicken, Jason E and Zingg, David W},
  journal={Journal of Computational and Applied Mathematics},
  volume={237},
  number={1},
  pages={111--125},
  year={2013},
  publisher={Elsevier}
}

@inproceedings{spiegel2015survey,
  title={A survey of the isentropic {E}uler vortex problem using high-order methods},
  author={Spiegel, Seth C and Huynh, HT and DeBonis, James R},
  booktitle={22nd AIAA computational fluid dynamics conference},
  pages={2444},
  year={2015}
}

@article{rydin2018high,
  title={High-fidelity sound propagation in a varying {3D} atmosphere},
  author={Rydin, Ylva and Mattsson, Ken and Werpers, Jonatan},
  journal={Journal of Scientific Computing},
  volume={77},
  number={2},
  pages={1278--1302},
  year={2018},
  publisher={Springer}
}

@article{duan2025active,
  title={Active flux methods for hyperbolic conservation laws—flux vector splitting and bound-preservation},
  author={Duan, Junming and Barsukow, Wasilij and Klingenberg, Christian},
  journal={SIAM Journal on Scientific Computing},
  volume={47},
  number={2},
  pages={A811--A837},
  year={2025},
  publisher={SIAM}
}

\appendix

\revThree{
\section{Boundary optimized, diagonal-norm upwind SBP operators}\label{AppA}
}

\revThree{
The boundary optimized diagonal-norm upwind SBP operators are available in the Julia library SummationByPartsOperators.jl~\cite{ranocha2021sbp} 
The available \texttt{accuracy\_order} of the operators are between $2$ to $9$, inclusive.
For example, one can construct a sixth order boundary-optimized upwind SBP operator containing the matrix pair $D_+$ and $D_-$, the norm matrix $H$, and the grid point vector $\mathbf{x}$ via}\\

\begin{verbatim}
D_upwind = upwind_operators(MattssonNiemeläWinters2026,
                            derivative_order = 1,
                            accuracy_order = 6,
                            xmin = 0.0, xmax = 1.0,
                            N = 13)
\end{verbatim}
\revThree{
where one can adjust the domain limits \texttt{xmin}, \texttt{xmax} or the number of grid points \texttt{N} as desired.
There is a minimum number of points required to properly define the boundary closures depending on the order of the operator.
For instance, the sixth order operator requires $N > 12$ whereas the ninth order operator must have $N > 16$.
For more details, see} \url{https://ranocha.de/SummationByPartsOperators.jl/stable/introduction/#intro-upwind-operators}.

\revThree{
Additionally, they are available as MATLAB function files at} \url{https://github.com/andrewwinters5000/optimized_upwind_sbp_operators}. 
\revThree{The MATLAB functions return the differentiation matrices $D_+$ and $D_-$, the norm matrix $H$, the grid point vector $\mathbf{x}$, and the interior grid spacing $h$. Input to the functions are the number of grid points $m$ and the width $L$ of the domain.}

\end{document}